\documentclass[2p]{article}
\usepackage{amssymb,amsmath,amsthm}
\usepackage{indentfirst}
\usepackage{exscale}
\usepackage{relsize}
\usepackage[numbers,sort&compress]{natbib}

\usepackage{geometry}
\usepackage{color}
\newcommand{\R}{\mathbb{R}}
\newtheorem{theorem}{Theorem}[section]

\theoremstyle{definition}

\newtheorem{lemma}[theorem]{Lemma}

\newtheorem{corollary}[theorem]{Corollary}

\allowdisplaybreaks
\theoremstyle{remark}
\newtheorem{remark}[theorem]{Remark}
\numberwithin{equation}{section}

\begin{document}
\title{\Large\bf{ Multiplicity result on a class of  nonhomogeneous quasilinear  elliptic system with small perturbations in $\mathbb{R}^N$}
 }
\date{}
\author {\ Xingyong Zhang$^{1,2}$, \ Wanting Qi$^{1}$\footnote{Corresponding author, E-mail address: qiwanting1@163.com}\\
      {\footnotesize $^{1}$Faculty of Science, Kunming University of Science and Technology, Kunming, Yunnan, 650500, P.R. China.}\\
      {\footnotesize $^{2}$Research Center for Mathematics and Interdisciplinary Sciences, Kunming University of Science and Technology,}\\
 {\footnotesize Kunming, Yunnan, 650500, P.R. China.}\\
 }
 \date{}
 \maketitle

 \begin{center}
 \begin{minipage}{15cm}

\small  {\bf Abstract:}
We investigate a class of quasilinear  elliptic system involving a nonhomogeneous differential operator which is introduced by C. A. Stuart [Milan J. Math. 79 (2011), 327-341] and depends on not only $\nabla u$ but also $u$. We show that the existence of multiple small solutions when the nonlinear term $F(x,u,v)$ satisfies locally sublinear and symmetric conditions and the perturbation is any continuous function with a small coefficient and no any growth hypothesis.
Our technical approach is mainly based on a variant of Clark's theorem without the global symmetric condition.
We develop the Moser's iteration technique to this quasi-linear elliptic system with nonhomogeneous differential operators and obtain that the relationship between $\|u\|_{\infty}$, $\|v\|_{\infty}$ and $\|u\|_{2^{\ast}}$, $\|v\|_{2^{\ast}}$.
We overcome some difficulties which are caused by the nonhomogeneity of the differential operator and the lack of compactness of the Sobolev embedding.
\end{minipage}
 \end{center}
  \allowdisplaybreaks

 \vskip2mm
\section{Introduction}\label{section 1}
\noindent
In this paper, we are concerned with the following non-homogeneous elliptic system with perturbation:
\begin{equation}\label{eq1}
 \left\{
  \begin{array}{ll}
 -\mbox{div}\left\{\phi_{1} \left(\frac{V_{1}(x)u^{2}+|\nabla u|^{2}}{2}\right)\nabla u\right\}+\phi_{1} \left(\frac{V_{1}(x)u^{2}+|\nabla u|^{2}}{2}\right)V_{1}(x)u
 =F_u(x,u,v)+\varepsilon k(x)G_{u}(u,v),\;\; x\in \mathbb{R}^{N},
 \\
 -\mbox{div}\left\{\phi_{2} \left(\frac{V_{2}(x)v^{2}+|\nabla v|^{2}}{2}\right)\nabla v\right\}+\phi_{2} \left(\frac{V_{2}(x)v^{2}+|\nabla v|^{2}}{2}\right)V_{2}(x)v
 =
 F_v(x,u,v)+\varepsilon k(x)G_{v}(u,v), \;\; x\in \mathbb{R}^{N},
 \\
 u\in H^{1}(\mathbb{R}^{N}),\;\;v\in H^{1}(\mathbb{R}^{N}),
  \end{array}
 \right.
 \end{equation}
where $N>2$ is an integer, $\phi_{i}:[0,+\infty)\rightarrow \mathbb{R}$ $(i=1,2)$ are two continuous functions which satisfy the following conditions:
\begin{itemize}
\item[$(\Upsilon_{1})$] there exist two constants $0< \rho_0\leq \rho_1$ such that
 $ 0<\rho_0\leq \phi_{i}(s)\leq \rho_1$ for all $s\in [0,+\infty)$;
\item[$(\Upsilon_{2})$] let $\hbar_{i}(s):=\Phi_{i}(s^{2})$ where $\Phi_{i}(s):=\int_{0}^{s}\phi_{i}(\varsigma)d\varsigma$, there exists $l_{i}>0$ such that
$\hbar_{i}(t)\geq \hbar_{i}(s)+\hbar_{i}'(s)(t-s)+l_{i}(t-s)^{2}$
for all $t,s\geq0$;
\item[$(\Upsilon_{3})$] $\phi_{i}(s)s\geq \Phi_{i}(s)$ for all $s\in [0,+\infty)$.
\end{itemize}
\par
Quasilinear elliptic system \eqref{eq1} contains the following class of quasi-linear elliptic scalar  equation:
\begin{equation}\label{eq2}
\left\{
\begin{array}{ll}
-\mbox{div}\left\{\phi\left(\frac{V(x)u^{2}+|\nabla u|^{2}}{2}\right)\nabla u\right\}
+\phi \left(\frac{V(x)u^{2}+|\nabla u|^{2}}{2}\right)V(x)u
=f(x,u)+\varepsilon k(x)g(u), &x\in \mathbb{R}^N,\\
u \in H^{1}(\mathbb{R}^{N}),
\end{array}
\right.
\end{equation}
which corresponds to the special situation $\phi_2=\phi_1=:\phi$, $v=u$, $V_{2}(x)=V_{1}(x):=V(x)$, $F(x,u,v)=F(x,v,u)$,  $G(u,v)=G(v,u)$, $F_{u}(x,u,v)=f(x,u)$ and $G_{u}(u,v)=g(u)$ in (\ref{eq1}).
The problems like \eqref{eq1} or \eqref{eq2} usually appear in the study of nonlinear optics model which describes the propagation of self-trapped beam in a cylindrical optical fiber made from a self-focusing dielectric material.
For details, we refer the readers to \cite{Landau1984,Mihalache1989,Svelto1974} for the essential physical background, and to \cite{Chen1991,Snyder1991,Stuart1996,Stuart1997,Stuart2010,Stuart2001} for the procedure to study the propagation of self-trapped beam in nonlinear optics.
\par
When $\varepsilon=0$, there have been some results about the existence and multiplicity of
solutions for system \eqref{eq1} or equations like  \eqref{eq2} on a bounded domain $\Omega\subset \R^N$,
for example, \cite{Stuart2011,Jeanjean2022,Qi-Zhang2023,Qi-Zhang2024,zhangxingyong2023}.
In \cite{Stuart2011}, Stuart considered problem \eqref{eq2} with $\varepsilon=0$ and $f(x,u)=\lambda u+h$ under the Dirichlet boundary condition $ u(x)=0$, $ x\in \partial\Omega$, where $N\geq1$, $\lambda\in\mathbb{R}$,
$h\in L^{2}(\Omega)$ and $h\geq0$ a.e. in $\Omega$.
Via an improved mountain pass theorem and the definition of localizing the Palais-Smale sequence,
they obtained
the existence of two non-negative weak solutions for problem (\ref{eq2}): one is a local minimum of the corresponding variational functional and the other is mountain-pass solution.
In \cite{Jeanjean2022}, Jeanjean-R\u{a}dulescu considered
problem \eqref{eq2}
with $\varepsilon=0$ and $f(x,u)=f(u)+h$, which extend nonlinearities in \cite{Stuart2011} directly to the case of the nonlinear growth reaction term.
To be specific, $h$ is non-negative and $f$ is a given continuous function which has either a sublinear decay or a linear growth at infinity.
In the sublinear case, $\phi$ was supposed to satisfy some reasonable conditions, and the existence of non-negative solutions for problem (\ref{eq2}) was obtained by using a minimization procedure.
In the linear growth case, based on the mountain pass theorem and the Ekeland's variational theorem,
they proved problem (\ref{eq2}) has at least one or two non-negative solutions
under stronger assumptions on $\phi$.
In \cite{Qi-Zhang2023}, we considered the multiplicity of solutions for problem (\ref{eq2}) with concave-convex nonlinearities and sign-changing weight functions. By means of the Nehari manifold and doing a fine analysis associated on the fibering map, we obtained that problem (\ref{eq2}) admits at least one positive energy solution  and negative energy solution  which is also the ground state solution of problem (\ref{eq2}).
In \cite{Qi-Zhang2024}, we further considered the case of system (\ref{eq1}) with concave-convex and nonhomogeneous terms and the existence of two nontrivial solutions  was obtained by the mountain pass theorem and Ekeland's variational principle.
Recently, in \cite{zhangxingyong2023}, our first author and Yu considered the existence of ground state sign-changing solutions for problem (\ref{eq2}) with $\varepsilon=0$.
Combining a non-Nehari manifold method and Miranda theorem, they showed that
the problem has one ground state sign-changing solution with two precise nodal domains when the nonlinearity is a continuous function, which has a superlinear growth at  infinity and satisfies some reasonable conditions.
Moreover, they also obtained that the energy of the ground state sign-changing solutions is larger than twice of the energy of the ground state solutions.
For quasilinear problem \eqref{eq2} on the entire space $\mathbb{R}^N$,
as a result of the lack of compactness of the Sobolev embedding,
to the best of our knowledge,
there is only one paper to consider the existence and multiplicity of solutions for problem (\ref{eq2}) (see \cite{Pomponio2021}).
In \cite{Pomponio2021}, Pomponio-Watanabe studied
problem (\ref{eq2}) with $\varepsilon=0$ and general nonlinear terms of Berestycki-Lions' type.
By using the mountain pass theorem and a technique of adding one dimension for space $\mathbb{R}^{N}$,
they proved that problem (\ref{eq2}) possesses at least a non-trivial weak solution
when the nonlinearity is a continuous and odd function, which has at least a linear growth near the origin.
They also obtained the existence of a radial ground state solution and a ground state solution by establishing the regularity of solutions and the Pohozaev identity.
\par
If $\varepsilon$ is allowed to be non-zero, we turn our gaze to a special example of equation \eqref{eq2} with $\phi(t)\equiv1$ (or more generally when $\phi(t)$ is a positive constant) and $V(x)\equiv0$,
\begin{equation}\label{eq3}
 \left\{
  \begin{array}{ll}
 -\Delta u=f(x,u)+\varepsilon g(x,u), &x\in \Omega,\\
 u=0, &x\in\partial\Omega,
 \end{array}
 \right.
\end{equation}
where $\Omega$ is a bounded domain in $\mathbb{R}^N$ with smooth boundary $\partial \Omega$ and $\varepsilon$ is a small parameter.
When the the nonlinear term $f(x,u)$ and $g(x,u)$ are sublinear and odd only near $u=0$,
the variational functional is not well defined and continuously differentiable on the working space,
which cause that the classical Clark's Theorem (\cite{Rabinowitz1986}, Theorem 9.1) can not work well for equation (\ref{eq2}).
Such problem could be dealt with well by some variant of the Clark Theorem, for example,
\cite[Theorem A. and Theorem 1.1]{Liu-Wang2015}.
If $g(x,u)$ is allowed to be any continuous function near $u=0$ without any other hypothesis,
the symmetry of equation is broken by the perturbation term $\varepsilon g$, which cause the above method is no longer valid.
This challenging case was considered by Kajikiya in \cite{Kajikiya2013} and the existence of at least $k$ distinct solutions of equation \eqref{eq3} for any given $k\in \mathbb N$ was obtained by developing a new version of Clark's theorem and using cut-off technique.
Subsequently, the idea in \cite{Kajikiya2013} has been applied to some different problems, for examples,
semilinear Schr\"{o}dinger equations \cite{Zhang-Zhang2024},
$p$-Kirchhoff type equations \cite{Liu-Wang-Zhao2016},
quasilinear elliptic equations or system \cite{Huang2022,Liu-Zhang-Wang2024} and so on.
\par
Based on the works in \cite{Jeanjean2022,Pomponio2021,Kajikiya2013,Liu-Zhang-Wang2024},
the motivation of our work is to consider whether the variant of Clark's theorem without the global symmetric condition due to Kajikiya in \cite{Kajikiya2013} can be applied to the quasilinear system \eqref{eq1} which is defined on the whole-space $\mathbb{R}^N$.
We overcome some difficulties which are caused by the nonhomogeneity of the differential operator, the coupling relationship of $u$ and $v$ and the lack of compactness of the Sobolev embedding.
Especially, we develop the Moser's iteration technique to this quasi-linear elliptic system with nonhomogeneous differential operators.
\par
There are three key steps in our proofs, which is on the basis of the idea in \cite{Kajikiya2013}:
(1) the construction of modified problem corresponding to system (\ref{eq1}).
Since we suppose that $F$ has a sublinear decay near the origin and is locally even with respect to $(u,v)$
and the perturbation term $G$ is any continuous function with a small coefficient $\varepsilon$.
the functionals $\int_{\R^{N}} F(x,u,v)dx$ and $\int_{\mathbb{R}^{N}}k(x)G(u,v)dx$ are not well defined in the workspace.
With the help of the cut-off function, we could extend $F$ and $G$ to two proper functions $\bar{F}$ and $\bar{G}$, respectively.
We replace $F$ and $G$ in system (\ref{eq1}) with  $\bar{F}$ and $\bar{G}$ to get a new problem, which is called the modified problem corresponding to system (\ref{eq1}).
By the original conditions for $F$ and $G$ and the definition of the cut-off function, we could obtain some global properties of $\bar{F}$ and $\bar{G}$, which make the variational functional corresponding to the modified problem  well defined and of class $C^{1}$.
(2) The solutions of modified problem. We will use Lemma \ref{lemma2.2} to obtain a multiplicity result of nontrivial weak solutions for modified problem.
Specifically, we shall prove that the variational functional corresponding to the modified problem satisfies all conditions of Lemma \ref{lemma2.2}
and obtain that for any $k\in\mathbb{N}$, there exists an constant $\varepsilon(k)>0$ such that if $|\varepsilon|\leq \varepsilon(k)$, the modified problem possesses at least $k$ distinct solutions.
(3) The solutions of the original problem \eqref{eq1}.
It is worth noting that these $k$ solutions are not necessarily solutions to the original problem \eqref{eq1}.
According to the definition of the cut off function, once we can control the infinite norm of these solutions within an appropriate range, these solutions are natural solutions of the original problem \eqref{eq1}.
Therefore, the key lies in the estimation of the infinite norm of solutions.
By the Moser's iteration technique and the continuity of the Sobolev embedding,
we can  establish  the relationship between the infinite norm of solutions and the norm for the workspace of solutions.
Then the estimation of the infinite norm of solutions can be translated into the estimation of the norm of the workspace of solutions (see Remark \ref{remark3.9.1}).
A crucial lemma (Lemma \ref{lemma3.7}) implies that the estimation of the norm of the workspace of solutions can be achieved by estimating the critical values (see Remark \ref{remark3.7.1}).
Fortunately, the critical values we obtained are negative and has an arbitrarily small lower boundary (see (\ref{z.1})).
We can take the appropriate lower bound such that these $k$ solutions are also solutions of the original problem \eqref{eq1}.
\par
Our results develop those results in  \cite{Kajikiya2013,Jeanjean2022,Pomponio2021,zhangxingyong2023,Liu-Zhang-Wang2024} from the following aspects:
\begin{itemize}
\item[(\uppercase\expandafter{\romannumeral1})]
Different from  \cite{Kajikiya2013},
we work in the entire space $\mathbb{R}^N$ rather than on a bounded domain $\Omega\subset\mathbb{R}^N$
and we focus on the more general quasilinear system  \eqref{eq1} with non-homogeneous differential operator instead of equation  \eqref{eq3} with homogeneous differential operator.
As a result of the coupling relationship of $u$ and $v$ and the non-homogeneity of   our operator,
proofs in the present paper become more difficult and complex than those in \cite{Kajikiya2013}.
Especially, such difficult and complexity can be embodied in two ways:
(i) establishment of the relationship between $\|u\|_{\infty}$, $\|v\|_{\infty}$ and $\|u\|_{2^{\ast}}$, $\|v\|_{2^{\ast}}$ (Lemma \ref{lemma3.8} below);
(ii) the proofs of $(PS)_{c}$-condition holds for functional $J_{\varepsilon}$ (Lemma \ref{lemma3.5}, Lemma \ref{lemma3.3} and Lemma \ref{lemma3.5} below).
To solve the first obstacle,  we develop the Moser's iteration technique to this quasi-linear elliptic system \eqref{eq1} and use some properties of our operator sufficiently to scale some inequalities carefully.
We overcome the second difficulty with the help of the idea of Stuart in \cite{Stuart2011} and the condition $(\Upsilon_{2})$.
\item[(\uppercase\expandafter{\romannumeral2})]
Our conditions are different from those in \cite{Jeanjean2022,Pomponio2021,zhangxingyong2023} if $\varepsilon=0$ because all of our assumptions on $F$ are local,  just near the origin, but those  assumptions in \cite{Jeanjean2022,Pomponio2021,zhangxingyong2023} are global.
It is worth noting that the existence of infinitely many solutions is not considered in \cite{Jeanjean2022,Pomponio2021,zhangxingyong2023} and this gap is filled in the present paper.
\item[(\uppercase\expandafter{\romannumeral3})]
For quasilinear equation \eqref{eq2}, on the whole space $\mathbb{R}^N$, the main difficulty is the lack of compactness of the Sobolev embedding,
which is crucial to ensure $(PS)_{c}$-condition holds for variational functional.
A useful way to overcome this difficulty is to reconstruct the compactness embedding theorem, which can be done by choosing the radially symmetric function space as the working space ( for example, see \cite{Pomponio2021}).
Being different from the way in \cite{Pomponio2021} and motivated by \cite{Liu-Zhang-Wang2024},
in this paper, we choose a subspace depending on $V$ as the working space where the function $V(x)^{-1}$ belongs to $L^{\alpha}(\mathbb{R}^{N})\cap L^{1}(\mathbb{R}^{N})$ ($\alpha$ will be given in the assumption $(\Lambda_{1})$ below), which is not a radially symmetric function space.
\item[(\uppercase\expandafter{\romannumeral4})]
Our results are different from those in \cite{Liu-Zhang-Wang2024}.
In \cite{Liu-Zhang-Wang2024}, Liu et al extended the result in \cite{Kajikiya2013} to the quasilinear elliptic system with a class of inhomogeneous differential operator called as the $\Phi$-Laplacian operator.
This operator is defined by the conditions $(\phi_1)$--$(\phi_3)$ in \cite{Liu-Zhang-Wang2024} and depends only on $\nabla u$ and usually is studied in the Orlicz-Sobolev space.
Our hypothesis $(\Upsilon_1)$ means that the function $\phi(t)$ lies between two positive constants, which is incompatible with conditions $(\phi_1)$--$(\phi_3)$ in \cite{Liu-Zhang-Wang2024}.
Therefore, our problem different from those in  \cite{Liu-Zhang-Wang2024} and our results cover some situations (see Remark \ref{remark1.1.1} and Section \ref{section 4} ) which does not satisfy conditions $(\phi_1)$--$(\phi_3)$ in \cite{Liu-Zhang-Wang2024}.
Although there are some similarities with \cite{Liu-Zhang-Wang2024} in terms of ideas and our proofs only needs to be conducted in the Hilbert space rather than the Orlicz-Sobolev space thanks to the assumptions $(\Upsilon_1)$--$(\Upsilon_3)$, our operator depends not only on $\nabla u$ but also on $u$, which causes that our proofs are different from those in \cite{Liu-Zhang-Wang2024}, especially, in the proofs that $(PS)_{c}$-condition holds for the variational functional.
\end{itemize}

\par
\noindent
\begin{remark}\label{remark1.1.0}
Here we make some supplement to conditions $(\Upsilon_{1})$--$(\Upsilon_{3})$.
\begin{itemize}
\item[(\uppercase\expandafter{\romannumeral1})]
The condition $(\Upsilon_{2})$ was introduced by Stuart in \cite{Stuart2011} to ensure that the following two aspects of the facts are true: (1) the convergence of a special class of bounded Cerami sequence $\{u_{n}\}$; (2) the uniform ellipticity of problem.
For the first aspect, the reflexivity and separability of the workspace ensure that the bounded Cerami sequence $\{u_{n}\}$ has weakly convergent subsequence.
Due to the completeness of the workspace, the proof of that this weak convergence is strong convergence can actually be converted into proof that Cerami sequence $\{u_{n}\}$ is a Cauchy sequence in the workspace, i.e., for any $n,m>0$, $\|u_{n}-u_{m}\|$ can be arbitrarily small.
Stuart established an inequality about $\|u_{n}-u_{m}\|$ by using the condition  $(\Upsilon_{2})$.
Specifically, the left-hand side of inequality is $\|u_{n}-u_{m}\|$ and the right-hand side of the inequality is related to the variational functional.
Then by the definition of Cerami sequence, linear growth behavior of the reaction term $\lambda u$ and the compactness of the Sobolev embedding, an arbitrarily small result of $\|u_{n}-u_{m}\|$ was obtained.
Inspired by Stuart's ideas, in this paper, we use the condition $(\Upsilon_{2})$ to prove the convergence of a class of bounded $(PS)_{c}$ sequence.
Since we investigate systems with general nonlinear terms on unbounded domains
rather than equations with linear growth term $\lambda u$ on bounded domains,
the coupling relationship of $u$ and $v$ and the lack of compactness of the Sobolev embedding make proofs in the present paper become more complex than those in \cite{Stuart2011}.
\item[(\uppercase\expandafter{\romannumeral2})]
The condition $(\Upsilon_{1})$ was introduced by Jeanjean and R\u{a}dulescu in \cite{Jeanjean2022}, which was inspired by condition $(g1)$ in \cite{Stuart2011}:
\begin{itemize}
\item[$(g1)$] $\phi$ is non-increasing on $[0,+\infty)$, $\phi(\infty)=\lim_{t\rightarrow\infty}\phi(t)>0$ and $\phi(0)=\lim_{t\rightarrow0}\phi(t)>0$.
\end{itemize}
Clearly, $(\Upsilon_{1})$ is weaker than the condition $(g1)$.
$(g1)$ implies that the function $\phi$ is bounded between two positive constants $\phi(\infty)$ and $\phi(0)$, which was used to scale some inequalities  in \cite{Stuart2011}.
Moreover, $(g1)$ requires $\Phi(t)$ to be a concave function of $t$.
The concavity of $\Phi(t)$ was used to prove the existence of a local minimum, the boundedness of a special class of bounded Cerami sequence $\{u_{n}\}$ and so on in \cite{Stuart2011}.
Being different from \cite{Stuart2011}, in this paper, we are concerned with the existence of infinitely many solutions for system (\ref{eq1}) rather than the existence of a local minimum.
The main tool in this paper is a variant of Clark's theorem without the global symmetric condition instead of the Ekeland's variational theorem  used in \cite{Stuart2011},
and in this paper, the concavity of $\Phi(t)$ will be no longer needed. Moreover,
we just need to require $(\Upsilon_{1})$ instead of $(g1)$.
\item[(\uppercase\expandafter{\romannumeral3})]
The condition $(\Upsilon_{3})$ was introduced for the first time in this paper, which is inspired by the following condition $(\phi_{3})$ in \cite{Pomponio2021,Qi-Zhang2024}:
\begin{itemize}
\item[$(\phi_3)$] $\Phi_{i}(s)\geq \phi_{i}(s)s$ for all $s\in [0,+\infty)$, $(i=1,2)$.
\end{itemize}
Obviously, these two conditions are different.
Consider $\phi(s)=-(1+s)^{-\alpha}+A$ for $s\in [0,+\infty)$ with $A>0$ and $\alpha>0$.
Then $(\Upsilon_{3})$ is satisfied and $(\phi_3)$ fails for any  $\alpha\in (0,1]$.
$(\Upsilon_{3})$ always fails and $(\phi_3)$ is satisfied for any $\alpha\in (1,+\infty)$.
\end{itemize}
\end{remark}

\par
To state our results, we introduce the following  assumptions concerning $F$, $V_{i}$, $k$ and $G$:
\begin{itemize}
\item[$(F_0)$]
$F: \mathbb{R}^N\times [-\delta, \delta] \times [-\delta, \delta]\rightarrow \mathbb{R}$ is a $C^1$ function for some $\delta>0$,  such that
$F(x,0,0)=0$ for all $x\in \mathbb{R}^N$;
\item[$(F_{1})$] there exist two  constants $C_1, C_2>0$ such that
\begin{equation*}
 \left\{
  \begin{array}{l}
 |F_t(x,t,s)|\leq C_1\left(|t|^{k_{1}r-1}+|s|^{\frac{k_{2}(k_{1}r-1)}{k_{1}}}\right),\\
 |F_s(x,t,s)|\leq C_2\left(|t|^{\frac{k_{1}(k_{2}r-1)}{k_{2}}}+|s|^{k_{2}r-1}\right),
    \end{array}
 \right.
 \end{equation*}
for all $|(t,s)|\leq\delta$  and $x\in  \mathbb{R}^{N}$, where $2\leq k_{i}<2^{\ast}$($i=1,2$), $r$ is a constant with $0<r<1$ such that $1<k_{i}r<\min\{2,\frac{2^{\ast}}{2}\}$($i=1,2$) and $k_{1}+k_{2}-k_{1}k_{2}r>0$.
\item[$(F_{2})$] there exists a constant $\beta\in(0,2)$ such that
$$
tF_{t}(x,t,s)+ sF_{s}(x,t,s)-\beta F(x,t,s)
<
\frac{\rho_{0}(2-\beta)}{4}(V_{1}(x)|t|^{2}+V_{2}(x)|s|^{2})
$$
for all $|(t,s)|\leq\delta$  and $x\in  \mathbb{R}^{N}$;
\item[$(F_{3})$]
$$\lim_{|(t,s)|\rightarrow0}
\left(\inf_{x\in \mathbb{R}^{N}}\frac{F(x,t,s)}{|t|^{\beta}+|s|^{\beta}} \right)=+\infty;$$
\item[$(F_4)$]
 $F(x,-t,-s)=F(x,t,s)$
for all $|(t,s)|\leq\delta$  and $x\in  \mathbb{R}^{N}$;
\end{itemize}
\begin{itemize}
\item[$(K)$]
$k(x)\in L^{1}(\mathbb{R}^N)\cap L^{\infty}(\mathbb{R}^N)$;
\end{itemize}
\begin{itemize}
\item[$(G)$]
$G: [-\delta, \delta] \times [-\delta, \delta]\rightarrow \mathbb{R}$ is a $C^1$ function such that $G(0,0)=0$;
\end{itemize}
\begin{itemize}
\item[$(\Lambda_{0})$]
$V_{i}\in C(\mathbb{R}^{N},\mathbb{R})$
and there exists a constant $V_{0}$ such that
$\inf_{x \in \mathbb{R}^{N} }V_{i}(x)=V_{0}>0$,\; $i=1,2$;
\item[$(\Lambda_{1})$]
the function $[V_{i}(x)]^{-1}$ belongs to $L^{\frac{k_{i}r}{2-k_{i}r}}(\mathbb{R}^{N})\cap L^{1}(\mathbb{R}^{N})$, $i=1,2$.
\end{itemize}

\par
Our main result is as follows.
\par
\noindent
\begin{theorem}\label{theorem1.1}
Assume that $(\Upsilon_{1})$--$(\Upsilon_{3})$, $(F_0)$--$(F_4)$, $(\Lambda_{0})$, $(\Lambda_{1})$, $(K)$ and $(G)$ hold.
 Then for any $k\in\mathbb{N}$, there exists a constant $\varepsilon(k)>0$ such that if $|\varepsilon|\leq \varepsilon(k)$, system (\ref{eq1}) possesses at least $k$ distinct solutions whose $L^{\infty}$-norms are less than $\frac{\delta}{2}$.
\end{theorem}
\par
By Theorem \ref{theorem1.1}, it is easy to obtain the following corollary.
\par
\noindent
\begin{corollary}\label{corollary1.1}
 Assume that   $(\Upsilon_{1})$--$(\Upsilon_{3})$, $(F_0)$--$(F_4)$, $(\Lambda_{0})$ and $(\Lambda_{1})$ hold.
 Then system (\ref{eq1}) with $\varepsilon=0$ possesses infinitely many distinct solutions whose $L^{\infty}$-norms are less than $\frac{\delta}{2}$.
\end{corollary}

\par
\noindent
\begin{remark}\label{remark1.1.1} There are examples satisfying all our assumptions. For example,
consider $\phi(s)=(1+s)^{-\alpha}+A$ for $s\geq0$ with $A>0$ and $\alpha>0$.
Then $(\Upsilon_{1})$ is satisfied for any $A>0$ and $\alpha>0$.
For fixed $\alpha$, $(\Upsilon_{2})$ holds for large enough values of $A$ and $(\Upsilon_{3})$ holds if $\alpha>1$.  One can see the example about the nonlinear term $F$ in section 4 below.
\end{remark}
\par
The paper is structured as follows.
In section \ref{section 2} we review some useful facts for the involved functional spaces.
In section \ref{section 3} we provide the proofs of Theorem \ref{theorem1.1}.
In section \ref{section 4} we give an example that illustrate our results.
In section \ref{section 5} we give some results for the scalar equation.
\vskip2mm
 \noindent
{\bf Notation.} Let $N\geq 1$, $1\leq p\leq\infty$, $ L^{p}(\R^{N})$ denotes the usual Lebesgue space with norms
$\|u\|_{p}:=\left(\int_{\R^{N}}|u|^{p}dx\right)^{\frac{1}{p}},\;\;1\leq p<\infty$
and
$\|u\|_{\infty}:=\inf\{ C>0: |u(x)|\leq C \;\mbox{almost everywhere on $\R^{N}$}\}$,
the Sobolev space $W^{1,p}(\R^{N})=\{u\in L^{p}(\R^{N}):\frac{\partial u}{\partial x_{i}}\in L^{p}(\R^{N}),i=1,2,...,N\}$ with norm
$\|u\|_{1,p}:=\left(\int_{\R^{N}}(|\nabla u|^{p}+|u|^{p})dx\right)^{\frac{1}{p}},\;\;1\leq p<\infty$.
$W^{1,2}(\R^{N})$ is usually recorded as $H^{1}(\R^{N})$.
We denote the dual space of $X$ by $X^{\ast}$, the weak convergence by $\rightharpoonup$ and the strong convergence by $\rightarrow$.

\section{Preliminaries}\label{section 2}
This section focuses on a survey of concepts and results from Sobolev spaces that
will be used in the text  and introduce some results from variational methods.  For a deeper understanding of these concepts and results, we refer readers for more details to the books \cite{Struwe1996,Brezis2020,Adams2003,Rabinowitz1986}.
\par
If we assume that
$V_{i}\in C(\mathbb{R}^{N},\mathbb{R})$ which satisfies $\inf_{x \in \mathbb{R}^{N} }V_{i}(x)=V_{0}>0$,\; $i=1,2$,
then, in order to deal with the system (\ref{eq1}), we introduce the subspace $W_{i}$ of $H^{1}(\R^{N})$, which is defined by \begin{equation*}
 W_{i}=\left\{u\in H^{1}(\mathbb{R}^{N}):\int_{\mathbb{R}^{N}} V_{i}(x)u^{2}dx<\infty   \right\}
  \end{equation*}
with the norm
\begin{equation*}
\|u\|_{ W_{i}}=\left(\int_{\mathbb{R}^{N}} (|\nabla u|^{2}+V_{i}(x)u^{2})dx\right)^{\frac{1}{2}},\;\;i=1,2.
  \end{equation*}
It is easy to see that $(W_{i},\|\cdot\|_{W_{i}})$($i=1,2$) are two separable and reflexive Banach spaces (see \cite{Medeiros2008}).
\par
In what follows, we give the following embedding results, which will be used for several times in Section \ref{section 3}.
\noindent
\begin{remark}\label{remark2.1}
Under condition $(\Lambda_{0})$,  we have
$$
     W_{i}\hookrightarrow L^{p_i}(\mathbb{R}^{N}),~~i=1,2,
 $$
with  continuous embeddings if $2\leq p_i\leq 2^{\ast} (i=1,2)$ and  compact embeddings if $2\leq p_i< 2^{\ast} (i=1,2)$ and $V_{i}(i=1,2)$ satisfy condition $(\Lambda_{1})$.
As a result, there exists positive constant $S_{p_{i}}$ such that
\begin{eqnarray*}\label{3.1.9}
       \|u\|_{p_{i}}\leq S_{p_{i}}\|u\|_{W_{i}},\;\;\forall p_i \in [2,2^{\ast}],\;i=1,2
\end{eqnarray*}
and the embeddings
$$
W_{i}(B_R)\hookrightarrow L^{p_i}(B_R), ~i=1,2
$$
are compact, where $R>0$, $B_R(0):=\left\{x\in \mathbb{R}^N: |x|<R\right\}$ and $\|u\|_{p_{i}}=\left(\int_{\mathbb{R}^{N}}|u|^{p_{i}}dx\right)^{\frac{1}{p_{i}}}$.
These facts can be found in \cite{Medeiros2008,Shubin1999,Omana1992}.
\end{remark}
\par
Next, we give the following Gagliardo-Nirenberg-Sobolev inequality, which will be used in the proof of Lemma \ref{lemma3.8}.
\noindent
\begin{remark}\label{remark2.2}(For example, \cite[Theorem 9.9]{Brezis2020})
Let $1\leq p<N$, $p^{\ast}=\frac{Np}{N-p}$.
Then
$$
W^{1,p}(\R^{N})\hookrightarrow L^{p^{\ast}}(\R^{N}),
$$
and there exists a constant $D$, depending only on $p$ and $N$, such that
\begin{eqnarray}\label{r2.2}
       \|u\|_{p^{\ast}}\leq D\|\nabla u\|_{p},\;\;\forall u\in W^{1,p}(\R^{N}).
\end{eqnarray}
\end{remark}
\vskip2mm
\par
Finally, we recall a  version of Clark's theorem which was introduced in \cite{Kajikiya2013},
which will be used to prove our Theorem \ref{theorem1.1} in Section \ref{section 3}.
\par
\noindent
\begin{lemma}\cite{Kajikiya2013}\label{lemma2.2}
Let $W$ be an infinite dimensional Banach space. For any $\varepsilon\in[0,1]$,  $I_{\varepsilon}\in C(W, \mathbb{R})$. Suppose that  $I_{\varepsilon}(u)$ has a continuous partial derivative $I_{\varepsilon}'$ with respect to $u$ and  satisfies $(A_{1})-(A_{5})$ below.
\begin{itemize}
\item[$(A_{1})$:]
$\inf \{I_{\varepsilon}(u):\varepsilon \in [0,1], u\in W \}>-\infty$;
\item[$(A_{2})$:]
For $u\in W$, $|I_{\varepsilon}(u)-I_{0}(u)|\leq \psi(\varepsilon)$, where $\psi\in C([0,1],\mathbb{R})$ and $\psi(0)=0$, $I_{0}(u)=:I_{\varepsilon}(u)|_{\varepsilon=0}$;
\item[$(A_{3})$:]
$I_{\varepsilon}(u)$ satisfies the (PS)-conditions uniformly on $\varepsilon$, i.e. if a sequence $(\varepsilon_{k},u_{k})$ in $[0,1]\times W$ satisfies that
$\sup_{k}|I_{\varepsilon_{k}}(u_{k})|<\infty$ and $I_{\varepsilon_{k}}'(u_{k})$ converges to zero, then
$(\varepsilon_{k},u_{k})$ has a convergent subsequence;
\item[$(A_{4})$:]
$I_{0}(u)=I_{0}(-u)$ for $u\in W$ and $I_{0}(0)=0$;
\item[$(A_{5})$:]
For $u\in W\backslash\{0\}$, there exists a unique $\vartheta(u)>0$ such that $I_{0}(\vartheta u)<0$ if $0<|\vartheta|<\vartheta(u)$ and $I_{0}(\vartheta u)\geq 0$ if $|\vartheta|\geq \vartheta(u)$.
\end{itemize}
 Denote
$$
S^{k}:=\{x\in \mathbb{R}^{k+1}:|x|=1\},
$$
$$
A_{k}:=\{\alpha\in C(S^{k},W):\alpha \;\;\mbox{is odd}\},
$$
$$
d_{k}:=\inf_{\alpha\in A_{k}}\max_{x\in S^{k}}I_{0}(\alpha(x)).
$$
Let $k\in \mathbb{N}\backslash \{0\}$ satisfying $d_{k}<d_{k+1}$. Then there exist two constants $\varepsilon_{k+1}$, $c_{k+1}$ such that $0<\varepsilon_{k+1}\leq 1$, $d_{k+1}\leq c_{k+1}<-\psi(\varepsilon)$ for $\varepsilon \in [0,\varepsilon_{k+1}]$,  and for any $\varepsilon \in [0,\varepsilon_{k+1}]$, $I_{\varepsilon}(\cdot)$ has a critical value in the interval $[d_{k+1}-\psi(\varepsilon),c_{k+1}+\psi(\varepsilon)]$.
\end{lemma}
\par
\noindent
\begin{remark}\label{remark2.2}
Through checking the proof of Lemma \ref{lemma2.2}, it is not difficult to verify that Lemma \ref{lemma2.2} is also true if $(A_{3})$ is replaced by the following condition:
\begin{itemize}
\item[$(A_{3})'$:] For any given $c\in\R$,
$I_{\varepsilon}(u)$ satisfies  the $(PS)_{c}$-condition uniformly on $\varepsilon$, i.e. if a sequence $(\varepsilon_{k},u_{k})$ in $[0,1]\times W$ satisfies that
$I_{\varepsilon_{k}}(u_{k})$ converges to $c$ and $I_{\varepsilon_{k}}'(u_{k})$ converges to zero, then
$(\varepsilon_{k},u_{k})$ has a convergent subsequence.
\end{itemize}
\end{remark}

\vskip2mm
\section{Proofs}\label{section 3}
\setcounter{equation}{0}
Throughout this section, to apply Lemma \ref{lemma2.2}, we work in the  subspace $W:=W_{1}\times W_{2}$ of the space $H^{1}(\R^{N})\times H^{1}(\R^{N})$ with the norm
\begin{eqnarray}\label{3.1.1}
        \|(u,v)\|
  =    \|u\|_{W_{1}}+ \|v\|_{W_{2}}.
\end{eqnarray}
It is easy to see that  $(W,\|\cdot\|)$ is a separable and reflexive Banach space.
Then, we define the variational functional corresponding to system (\ref{eq1}) by
 \begin{eqnarray*}
 J_{\varepsilon}(u,v)
& :=&
        \int_{\mathbb{R}^{N}}\left\{\Phi_{1}\left(\frac{V_{1}(x)u^{2}+|\nabla u|^{2}}{2}\right)
        +\Phi_{2}\left(\frac{V_{2}(x)v^{2}+|\nabla v|^{2}}{2}\right)\right\}dx
\nonumber\\
 & &
        - \int_{\mathbb{R}^{N}}F(x,u,v)dx
        - \varepsilon\int_{\mathbb{R}^{N}}k(x)G(u,v)dx,\quad (u,v)\in W.
  \end{eqnarray*}
Conditions $(F_0)$, $(F_1)$ and $(G)$ imply the behaviors of  $F$ and $G$ are just near the origin. So the functionals $\int_{\R^{N}} F(x,u,v)dx$ and $\int_{\mathbb{R}^{N}}k(x)G(u,v)dx$ is not well defined in the Sobolev space $W$, respectively.
To deal with this problem, we extend $F$ and $G$ to two proper functions $\bar{F}$ and $\bar{G}$ by the cut-off technique, respectively,
which developed by Costa-Wang \cite{Costa-Wang2005} and has been applied to some different problems,
for example, \cite{Liu-Wang2015,Huang2022,Kajikiya2013,Liu-Zhang-Wang2024}.
\par
For some $\delta>0$, let $\tau\in C^1(\R ^{2},[0,1])$ as an even cut-off function defined by
 \begin{eqnarray}\label{3.0.1}
 \tau(t,s)= \begin{cases}
           1, \;\;\;  \text { if }\;\;|(t,s)| \leq \delta/2, \\
           0, \;\;\;  \text { if }\;\;|(t,s)|    >   \delta,
          \end{cases}
   \end{eqnarray}
which fulfills $t\tau'_{t}(t,s)\leq 0$ and $s\tau'_{s}(t,s)\leq 0$ for all $(t,s)\in\R ^{2}$ with $|(t,s)|\leq\delta$.
One can see some examples of $\tau(t,s)$ in \cite{Liu-Zhang-Wang2024}.
The modified functions $\bar{F}$ and $\bar{G}$ defined by $\bar{F}(x,t,s)=\tau(t,s)F(x,t,s)$ and $\bar{G}(t,s)=\tau(t,s)G(t,s)$, respectively.
Then, we have
$$
\bar{F}_{t}(x,t,s)=\tau_{t}(t,s)F(x,t,s)+\tau(t,s)F_{t}(x,t,s),\;\;
\bar{F}_{s}(x,t,s)=\tau_{s}(t,s)F(x,t,s)+\tau(t,s)F_{s}(x,t,s),
$$
$$
\bar{G}_{t}(t,s)=\tau_{t}(t,s)G(t,s)+\tau(t,s)G_{t}(t,s),\;\;
\bar{G}_{s}(t,s)=\tau_{s}(t,s)G(t,s)+\tau(t,s)G_{s}(t,s).
$$
It is easy to prove that $\bar{F}$ and  $\bar{G}$ satisfy the following conditions.
\par
\noindent
\begin{lemma}\label{lemma3.1}
Assume that $(F_0)$--$(F_4)$ and $(G)$ hold.
Then $\bar{F}$ and $\bar{G}$ are continuous functions and the following conditions hold:
\begin{itemize}
\item[$(F'_{0})$] $\bar{F}: \mathbb{R}^{N}\times \mathbb{R}\times \mathbb{R} \rightarrow \mathbb{R}$ is a $C^{1}$ function which satisfies $\bar{F}(x,0,0)=0$ for all $x\in \mathbb{R}^{N}$ and  $\bar{F}(x,t,s)\geq0$ for all  $(x, t,s)\in  \mathbb{R}^{N}\times \mathbb{R}\times \mathbb{R}$;
\item[$(F'_{1})$]
\begin{equation*}
 \left\{
  \begin{array}{l}
 |\bar{F}_t(x,t,s)|\leq C_1\left(|t|^{k_{1}r-1}+|s|^{\frac{k_{2}(k_{1}r-1)}{k_{1}}}\right),\\
 |\bar{F}_s(x,t,s)|\leq C_2\left(|t|^{\frac{k_{1}(k_{2}r-1)}{k_{2}}}+|s|^{k_{2}r-1}\right),
    \end{array}
 \right.
 \end{equation*}
for all  $(x, t,s)\in  \mathbb{R}^{N}\times \mathbb{R}\times \mathbb{R}$, where $2\leq k_{i}<2^{\ast}$($i=1,2$), $r$ is a constant with $0<r<1$ such that $1<k_{i}r<\min\{2,\frac{2^{\ast}}{2}\}$($i=1,2$) and $k_{1}+k_{2}-k_{1}k_{2}r>0$;
\item[$(F'_{2})$]
$$
t\bar{F}_{t}(x,t,s)+ s\bar{F}_{s}(x,t,s)-\beta \bar{F}(x,t,s)
<
\frac{\rho_{0}(2-\beta)}{2}(V_{1}(x)|t|^{2}+V_{2}(x)|s|^{2})
$$
for all $(x, t,s)\in  \mathbb{R}^{N}\times \mathbb{R}\times \mathbb{R}$, where $\beta<2$;
\item[$(F'_{3})$]
$$\lim_{|(t,s)|\rightarrow0}
\left(\inf_{x\in \mathbb{R}^{N}}\frac{\bar{F}(x,t,s)}{|t|^{\beta}+|s|^{\beta}} \right)=+\infty;$$
\item[$(F'_{4})$] $\bar{F}(x,-t,-s)=\bar{F}(x,t,s)$ for all $(x, t,s)\in  \mathbb{R}^{N}\times \mathbb{R}\times \mathbb{R}$;
\item[$(F'_{5})$] $\bar{F}(x,t,s)=\bar{F}_{t}(x,t,s)=\bar{F}_{s}(x,t,s)=\bar{G}(t,s)=\bar{G}_{t}(t,s)=\bar{G}_{s}(t,s)=0$ for all $x\in \mathbb{R}^{N}$ and $(t,s)\in \R^{2}$ with $|(t,s)|>\delta$.
\end{itemize}
\end{lemma}
\par
\noindent
\begin{remark}\label{remark3.4.0}
From condition $(F_3)$, we may assume that
 $$
 F(x,t,s)>0 \;\; \mbox{for all $0<|(t,s)|\leq \delta$ and $x\in \R^{N}$},
 $$
which is necessary to prove that $\bar{F}(x,t,s)\geq0$ for all  $(x, t,s)\in  \mathbb{R}^{N}\times \mathbb{R}\times \mathbb{R}$ and $(F'_{2})$ hold.
\end{remark}
\par
\noindent
\begin{remark}\label{remark3.4.1}
By the Young's inequality, $(F_0')$, $(F_1')$ and the fact
 $$\bar{F}(x,u,v)=\int_{0}^{u}\bar{F}_s(x,s,v)ds+\int_{0}^{v}\bar{F}_t(x,0,t)dt+\bar{F}(x,0,0), \quad \forall  (x,u,v)\in\mathbb{R}^{N}\times\mathbb{R}\times\mathbb{R},$$
it is easy to see that there exists a constant $C_3>0$ such that
\begin{eqnarray}\label{r3.4}
|\bar{F}(x,t,s)|\leq C_3(|t|^{k_{1}r}+|s|^{k_{2}r}),
\;\;\forall (x,t,s)\in \mathbb{R}^{N}\times \mathbb{R}\times \mathbb{R}.
\end{eqnarray}
\end{remark}
\par
\noindent
\begin{remark}\label{remark3.4.2}
From the assumptions about $k_{1}$, $k_{2}$ and $r$ in condition $(F'_{1})$, it is easy to verify that the following conclusions hold:
\begin{eqnarray}
& &
 \label{A3}
\left\{
\frac{k_1r-1}{k_1-1},
\frac{k_1-k_1r}{k_1-1},
\frac{k_2r-1}{k_2-1},
\frac{k_2-k_2r}{k_2-1}
\right\}
\subset (0,1),
\\
& &
 \label{A4}
\left\{
\frac{(2-k_{2}r)(k_{1}r-1)}{2k_{1}r},
\frac{(2-k_{1}r)(k_{2}r-1)}{2k_{2}r}
\right\}
\subset (0,1),
\\
& &
 \label{A2}
\left\{
\frac{k_{1}-k_{2}(k_{1}r-1)}{2k_{1}},
\frac{k_{2}(k_{1}r-1)}{2k_{1}},
\frac{k_{2}-k_{1}(k_{2}r-1)}{2k_{2}},
\frac{k_{1}(k_{2}r-1)}{2k_{2}}
\right\}
\subset (0,1),
\\
& &
 \label{A1}
\left\{
\frac{(k_{1}r-1)(2-k_{2}r)}{r(k_{1}-k_{2}(k_{1}r-1))},
\frac{2-k_{1}r}{r(k_{1}-k_{2}(k_{1}r-1))},
\frac{(k_{2}r-1)(2-k_{1}r)}{r(k_{2}-k_{1}(k_{2}r-1))},
\frac{2-k_{2}r}{r(k_{2}-k_{1}(k_{2}r-1))}
\right\}
\subset (0,1),
\\
& &\label{A5}
\left\{
\frac{k_{2}r(k_{1}-k_{2}(k_{1}r-1))}{2k_{1}(2-k_{2}r)},
\frac{k_{1}r(k_{1}-k_{2}(k_{1}r-1))}{2k_{1}(2-k_{1}r)},
\frac{k_{1}r(k_{2}-k_{1}(k_{2}r-1))}{2k_{2}(2-k_{1}r)},
\frac{k_{2}r(k_{2}-k_{1}(k_{2}r-1))}{2k_{2}(2-k_{2}r)}
\right\}
\subset (0,+\infty).\quad\quad
 \end{eqnarray}
\end{remark}
 \par
Consider the modified system of (\ref{eq1}) given by
\begin{equation}\label{mod1}
 \left\{
  \begin{array}{ll}
    -\mbox{div}\left\{\phi_{1} \left(\frac{V_{1}(x)u^{2}+|\nabla u|^{2}}{2}\right)\nabla u\right\}
    +\phi_{1} \left(\frac{V_{1}(x)u^{2}+|\nabla u|^{2}}{2}\right)V_{1}(x)u
 =
     \bar{F}_u(x,u,v)+\varepsilon k(x)\bar{G}_{u}(u,v),\;\; x\in \mathbb{R}^{N},
\\
    -\mbox{div}\left\{\phi_{2} \left(\frac{V_{2}(x)v^{2}+|\nabla v|^{2}}{2}\right)\nabla v\right\}
     +\phi_{2} \left(\frac{V_{2}(x)v^{2}+|\nabla v|^{2}}{2}\right)V_{2}(x)v
 =
    \bar{F}_v(x,u,v)+\varepsilon k(x)\bar{G}_{v}(u,v), \;\; x\in \mathbb{R}^{N},
\\
    u\in H^{1}(\mathbb{R}^{N}),\;\;v\in H^{1}(\mathbb{R}^{N}),
  \end{array}
 \right.
 \end{equation}
and define a functional $\bar{J}_{\varepsilon}$ on $W$ by
 \begin{eqnarray}\label{3.0.2}
        \bar{J}_{\varepsilon}(u,v)
& :=&
        \int_{\mathbb{R}^{N}}\left\{\Phi_{1}\left(\frac{V_{1}(x)u^{2}+|\nabla u|^{2}}{2}\right)
        +\Phi_{2}\left(\frac{V_{2}(x)v^{2}+|\nabla v|^{2}}{2}\right)\right\}dx
\nonumber\\
 & &
        - \int_{\mathbb{R}^{N}}\bar{F}(x,u,v)dx
        - \varepsilon\int_{\mathbb{R}^{N}}k(x)\bar{G}(u,v)dx,\quad (u,v)\in W.
 \end{eqnarray}
Under assumptions $(\Upsilon_{1})$, $(F_0')$, $(F_1')$, $(\Lambda_{0})$, $(\Lambda_{1})$, $(K)$ and $(G)$, by similar arguments as those in \cite[Lemma A.2]{wang2017}, we can prove that $\bar{J}_{\varepsilon}$ is well defined and of class $C^{1}(W,\mathbb{R})$ for all $\varepsilon \in \R$ with
 \begin{eqnarray}\label{3.0.3}
       \langle \bar{J}_{\varepsilon}'(u,v),(\varphi_{1},\varphi_{2})\rangle
 &= &
        \int_{\mathbb{R}^{N}}\phi_{1}
        \left(\frac{V_{1}(x)u^{2}+|\nabla u|^{2}}{2}\right)
        (V_{1}(x)u\varphi_{1}+\nabla u\cdot \nabla \varphi_{1})dx
\nonumber\\
 & &
       +\int_{\mathbb{R}^{N}}\phi_{2}
       \left(\frac{V_{2}(x)v^{2}+|\nabla v|^{2}}{2}\right)
       (V_{2}(x)v\varphi_{2}+\nabla v\cdot \nabla \varphi_{2} )dx
\nonumber\\
&    &
       -\int_{\mathbb{R}^{N}} \bar{F}_{t}(x,u,v)\varphi_{1}dx
       - \int_{\mathbb{R}^{N}}\bar{F}_{s}(x,u,v)\varphi_{2}dx
\nonumber\\
&    &
      -\varepsilon\int_{\mathbb{R}^{N}}k(x) \bar{G}_{t}(u,v)\varphi_{1}dx
      -\varepsilon\int_{\mathbb{R}^{N}}k(x) \bar{G}_{s}(u,v)\varphi_{2}dx
   \end{eqnarray}
for all $ (u,v), (\varphi_{1},\varphi_{2})\in W$.
Thus, the critical points of  $\bar{J}_{\varepsilon}$ in $W$ are weak solutions of system (\ref{mod1}).
\par
\noindent
\begin{lemma}\label{lemma3.4}
Assume that $(\Upsilon_{1})$, $(F_0')$, $(F_1')$, $(F_{5}')$, $(K)$, $(\Lambda_{0})$ and $(\Lambda_{1})$ hold.
Then the functional $\bar{J}_{\varepsilon}$ is coercive and bounded below on $W$ for any $\varepsilon\in [0,1]$.
\end{lemma}
\par
\noindent
{\bf Proof.}
For any given $(u,v)\in W$ and $\varepsilon\in [0,1]$,
it follows from  $(F'_{5})$ and $(K)$ that
\begin{eqnarray}\label{3.4.1.0}
      \varepsilon\int_{\mathbb{R}^{N}}|k(x)\bar{G}(u,v)|dx
&=&
        \varepsilon\int_{\{x\in \R^{N}: |(u(x),v(x))|\leq \delta\}}|k(x)||\bar{G}(u,v)|dx
      + \varepsilon\int_{\{x\in \R^{N}: |(u(x),v(x))|>\delta\}}|k(x)||\bar{G}(u,v)|dx
\nonumber\\
&=&
       \varepsilon\int_{\{x\in \R^{N}: |(u(x),v(x))|\leq \delta\}}|k(x)||\bar{G}(u,v)|dx
\nonumber\\
&\leq&
      \varepsilon\max_{|(t,s)|\leq \delta}|\bar{G}(t,s)|
     \int_{\{x\in \R^{N}: |(u(x),v(x))|\leq \delta\}}|k(x)|dx
\nonumber\\
&\leq&
      \varepsilon\max_{|(t,s)|\leq \delta}|\bar{G}(t,s)|\int_{\mathbb{R}^{N}}|k(x)|dx
\nonumber\\
&=&
     \varepsilon C_{4},
\end{eqnarray}
where $C_{4}=\max_{|(t,s)|\leq \delta}|\bar{G}(t,s)|\int_{\mathbb{R}^{N}}|k(x)|dx$.
Moreover, by (\ref{r3.4}), $1<k_{i}r<2$($i=1,2$),  H\"{o}lder's inequality and $(\Lambda_{1})$, we conclude that
\begin{eqnarray}\label{3.4.1.1}
& &
      \int_{\mathbb{R}^{N}}|\bar{F}(x,u,v)|dx
\nonumber\\
&\leq&
       C_3\int_{\mathbb{R}^{N}}(|u|^{k_{1}r}+|v|^{k_{2}r})dx
\nonumber\\
&=&
       C_3\int_{\mathbb{R}^{N}}
          (V_{1}^{-1}(x))^{\frac{k_{1}r}{2}}
          (V_{1}(x))^{\frac{k_{1}r}{2}}
          |u|^{k_{1}r}dx
      +C_3\int_{\mathbb{R}^{N}}
          (V_{2}^{-1}(x))^{\frac{k_{2}r}{2}}
          (V_{2}(x))^{\frac{k_{2}r}{2}}
          |v|^{k_{2}r}dx
\nonumber\\
&\leq&
        C_3
        \left(
        \int_{\mathbb{R}^{N}}
        (V_{1}^{-1}(x))^{\frac{k_{1}r}{2}\cdot\frac{2}{2-k_{1}r}}
         dx
        \right)^{\frac{2-k_{1}r}{2}}
        \left(
        \int_{\mathbb{R}^{N}}
        (V_{1}(x))^{\frac{k_{1}r}{2}\cdot\frac{2}{k_{1}r}}
        |u|^{k_{1}r\cdot\frac{2}{k_{1}r}}
        dx
        \right)^{\frac{k_{1}r}{2}}
\nonumber\\
&&
       +C_3
        \left(
        \int_{\mathbb{R}^{N}}
        (V_{2}^{-1}(x))^{\frac{k_{2}r}{2}\cdot\frac{2}{2-k_{2}r}}
         dx
        \right)^{\frac{2-k_{2}r}{2}}
        \left(
        \int_{\mathbb{R}^{N}}
        (V_{2}(x))^{\frac{k_{2}r}{2}\cdot\frac{2}{k_{2}r}}
        |v|^{k_{2}r\cdot\frac{2}{k_{2}r}}
        dx
        \right)^{\frac{k_{2}r}{2}}
\nonumber\\
&=&
       C_3
        \left(
        \int_{\mathbb{R}^{N}}
        (V_{1}^{-1}(x))^{\frac{k_{1}r}{2-k_{1}r}}
        dx
        \right)^{\frac{2-k_{1}r}{k_{1}r}\cdot\frac{k_{1}r}{2}}
        \left(
        \int_{\mathbb{R}^{N}}V_{1}(x)|u|^{2}dx
        \right)^{\frac{k_{1}r}{2}}
\nonumber\\
&&
       +C_3
        \left(
        \int_{\mathbb{R}^{N}}
        (V_{2}^{-1}(x))^{\frac{k_{2}r}{2-k_{2}r}}
         dx
        \right)^{\frac{2-k_{2}r}{k_{2}r}\cdot\frac{k_{2}r}{2}}
        \left(
        \int_{\mathbb{R}^{N}}V_{2}(x)|v|^{2}dx
        \right)^{\frac{k_{2}r}{2}}
\nonumber\\
&\leq&
        C_5
        \left(
        \int_{\mathbb{R}^{N}}(|\nabla u|^{2}+V_{1}(x)|u|^{2})dx
        \right)^{\frac{k_{1}r}{2}}
       +C_5
        \left(
        \int_{\mathbb{R}^{N}}(|\nabla v|^{2}+V_{2}(x)|v|^{2})dx
        \right)^{\frac{k_{2}r}{2}}
\nonumber\\
&=&
        C_5\|u\|_{W_{1}}^{k_{1}r}
       +C_5\|v\|_{W_{2}}^{k_{2}r},
\end{eqnarray}
where
\begin{eqnarray}\label{3.4.3}
C_5=C_3\max\left\{
\|V_{1}^{-1}\|_{\frac{k_{1}r}{2-k_{1}r}}^{\frac{k_{1}r}{2}},
\|V_{2}^{-1}\|_{\frac{k_{2}r}{2-k_{2}r}}^{\frac{k_{2}r}{2}}
\right\}.
 \end{eqnarray}
Thus, by (\ref{3.0.2}), $(\Upsilon_{1})$, $(\Lambda_{0})$, (\ref{3.4.1.0}) and (\ref{3.4.1.1}), we have
\begin{eqnarray}\label{3.4.2}
        \bar{J}_{\varepsilon}(u,v)
& =&
        \int_{\mathbb{R}^{N}}\left\{\Phi_{1}\left(\frac{V_{1}(x)u^{2}+|\nabla u|^{2}}{2}\right)
        +\Phi_{2}\left(\frac{V_{2}(x)v^{2}+|\nabla v|^{2}}{2}\right)\right\}dx
\nonumber\\
 & &
        - \int_{\mathbb{R}^{N}}\bar{F}(x,u,v)dx
        - \varepsilon\int_{\mathbb{R}^{N}}k(x)\bar{G}(u,v)dx
\nonumber\\
&\geq&
        \frac{\rho_{0}}{2}\|u\|_{W_{1}}^{2}
        +\frac{\rho_{0}}{2}\|v\|_{W_{2}}^{2}
        - \int_{\mathbb{R}^{N}}|\bar{F}(x,u,v)|dx
        - \varepsilon\int_{\mathbb{R}^{N}}|k(x)\bar{G}(u,v)|dx
\nonumber\\
&\geq&
        \frac{\rho_{0}}{2}\|u\|_{W_{1}}^{2}
        +\frac{\rho_{0}}{2}\|v\|_{W_{2}}^{2}
        -C_5\|u\|_{W_{1}}^{k_{1}r}
        -C_5\|v\|_{W_{2}}^{k_{2}r}
        -\varepsilon C_4,
 \end{eqnarray}
which implies that $\bar{J}_{\varepsilon}$ is coercive and bounded below on $W$ since $k_{i}r<2$ ($i=1,2$).
\qed

\vskip2mm
\par
Next, we prove that for any given $c\in\R$, the $(PS)_{c}$-condition holds for functional $\bar{J}_{\varepsilon}$. Until then,  we present the following lemmas which will be used in Lemma \ref{lemma3.5}.
\par
\noindent
\begin{lemma}\label{lemma3.2}
Assume that $(\Lambda_{0})$, $(\Lambda_{1})$ and $(F'_{1})$ hold.
If there exists a sequence  $(u_{n},v_{n}) \rightharpoonup (u,v)$ in $W$ as $n \rightarrow +\infty$, then
\begin{eqnarray*}
    \int_{\mathbb{R}^N}\left|\bar{F}_{t}(x,u_m,v_m)\right||u_n-u_m|dx\rightarrow 0
\;\mbox{and}\;
    \int_{\mathbb{R}^N}\left|\bar{F}_{s}(x,u_m,v_m)\right||v_n-v_m|dx\rightarrow 0
\end{eqnarray*}
as $n,m\rightarrow +\infty$.
\end{lemma}
\par
\noindent
{\bf Proof.}
By the condition $(\Lambda_{1})$, for any given constant $\epsilon>0$, there exists $R=R(\epsilon)>0$ such that
\begin{eqnarray}\label{3.2.2}
    \left(\int_{\mathbb{R}^{N}\backslash B_{R}(0)}
    \left(V_{i}(x)^{-1}\right)^{\frac{k_{i}r}{2-k_{i}r}}dx\right)^{\frac{2-k_{i}r }{k_{i}r}}
<   \epsilon,\;\;\forall\;i=1,2.
\end{eqnarray}
Let $\mathbb{A}:=\mathbb{R}^{N}\backslash B_{R}(0)$.
By
(\ref{A2}),  (\ref{A1}), (\ref{A5}) and Young's inequality, we obtain that
\begin{eqnarray}\label{3.2.3}
& &
       \left(
       \int_{\mathbb{A}}
       (V_{2}^{-1}(x))^{\frac{k_{2}(k_{1}r-1)}{k_{1}-k_{2}(k_{1}r-1)}}
       (V_{1}^{-1}(x))^{\frac{k_{1}}{k_{1}-k_{2}(k_{1}r-1)}}
       dx
       \right)^{\frac{k_{1}-k_{2}(k_{1}r-1)}{2k_{1}}}
\nonumber \\
&\leq&
      \Bigg\{
       \frac{(k_{1}r-1)(2-k_{2}r)}{r(k_{1}-k_{2}(k_{1}r-1))}
       \int_{\mathbb{A}}
       (V_{2}^{-1}(x))^{\frac{k_{2}(k_{1}r-1)}{k_{1}-k_{2}(k_{1}r-1)}
       \cdot \frac{r(k_{1}-k_{2}(k_{1}r-1))}{(k_{1}r-1)(2-k_{2}r)}}
       dx
\nonumber \\
&&
       +
       \frac{2-k_{1}r}{r(k_{1}-k_{2}(k_{1}r-1))}
       \int_{\mathbb{A}}
       (V_{1}^{-1}(x))^{\frac{k_{1}}{k_{1}-k_{2}(k_{1}r-1)}
       \cdot \frac{r(k_{1}-k_{2}(k_{1}r-1))}{2-k_{1}r}}
       dx
     \Bigg\}^{\frac{k_{1}-k_{2}(k_{1}r-1)}{2k_{1}}}
\nonumber \\
&=&
      \Bigg\{
       \frac{(k_{1}r-1)(2-k_{2}r)}{r(k_{1}-k_{2}(k_{1}r-1))}
       \int_{\mathbb{A}}
       (V_{2}^{-1}(x))^{\frac{k_{2}r}{2-k_{2}r}}
       dx
       +
       \frac{2-k_{1}r}{r(k_{1}-k_{2}(k_{1}r-1))}
       \int_{\mathbb{A}}
       (V_{1}^{-1}(x))^{\frac{k_{1}r}{2-k_{1}r}}
       dx
     \Bigg\}^{\frac{k_{1}-k_{2}(k_{1}r-1)}{2k_{1}}}
\nonumber \\
&\leq&
       \left(\frac{(k_{1}r-1)(2-k_{2}r)}{r(k_{1}-k_{2}(k_{1}r-1))}\right)
       ^{\frac{k_{1}-k_{2}(k_{1}r-1)}{2k_{1}}}
       \left(\int_{\mathbb{A}}(V_{2}^{-1}(x))^{\frac{k_{2}r}{2-k_{2}r}}dx\right)
       ^{\frac{k_{1}-k_{2}(k_{1}r-1)}{2k_{1}}}
\nonumber \\
&&
      +
      \left(\frac{2-k_{1}r}{r(k_{1}-k_{2}(k_{1}r-1))}\right)
      ^{\frac{k_{1}-k_{2}(k_{1}r-1)}{2k_{1}}}
       \left(\int_{\mathbb{A}}(V_{1}^{-1}(x))^{\frac{k_{1}r}{2-k_{1}r}}dx\right)
      ^{\frac{k_{1}-k_{2}(k_{1}r-1)}{2k_{1}}}
\nonumber \\
&=&
       \left(\frac{(k_{1}r-1)(2-k_{2}r)}{r(k_{1}-k_{2}(k_{1}r-1))}\right)
       ^{\frac{k_{1}-k_{2}(k_{1}r-1)}{2k_{1}}}
       \left(\int_{\mathbb{A}}(V_{2}^{-1}(x))^{\frac{k_{2}r}{2-k_{2}r}}dx\right)
       ^{\frac{2-k_{2}r}{k_{2}r}\cdot\frac{k_{2}r(k_{1}-k_{2}(k_{1}r-1))}{2k_{1}(2-k_{2}r)}}
\nonumber \\
&&
      +
      \left(\frac{2-k_{1}r}{r(k_{1}-k_{2}(k_{1}r-1))}\right)
      ^{\frac{k_{1}-k_{2}(k_{1}r-1)}{2k_{1}}}
       \left(\int_{\mathbb{A}}(V_{1}^{-1}(x))^{\frac{k_{1}r}{2-k_{1}r}}dx\right)
      ^{\frac{2-k_{1}r}{k_{1}r}\cdot\frac{k_{1}r(k_{1}-k_{2}(k_{1}r-1))}{2k_{1}(2-k_{1}r)}}
\nonumber \\
&<  &
     C
     \left(
     \epsilon^{\frac{k_{2}r(k_{1}-k_{2}(k_{1}r-1))}{2k_{1}(2-k_{2}r)}}
     +\epsilon^{\frac{k_{1}r(k_{1}-k_{2}(k_{1}r-1))}{2k_{1}(2-k_{1}r)}}
     \right).
 \end{eqnarray}
Thus, by  $(F'_{1})$, $1<k_{1}r<2$, (\ref{A2}), H\"{o}lder's inequality, (\ref{3.2.2}) and (\ref{3.2.3}), we get
\begin{eqnarray*}\label{3.2.4}
& &
      \int_{\mathbb{A}}|\bar{F}_{t}(x,u_{m},v_m)||u_{m}-u_{n}|dx
\nonumber\\
&\leq&
        C_1\int_{\mathbb{A}}
       \left(|u_{m}|^{k_{1}r-1}+|v_m|^{\frac{k_{2}(k_{1}r-1)}{k_{1}}}\right)
       (|u_{m}|+|u_n|)
       dx
\nonumber\\
&=&
     C_1\int_{\mathbb{A}}|u_{m}|^{k_{1}r}dx
     +C_1\int_{\mathbb{A}}|u_{m}|^{k_{1}r-1}|u_n|dx
     +C_1\int_{\mathbb{A}}|v_m|^{\frac{k_{2}(k_{1}r-1)}{k_{1}}}|u_{m}|dx
     +C_1\int_{\mathbb{A}}|v_m|^{\frac{k_{2}(k_{1}r-1)}{k_{1}}}|u_n|dx
\nonumber\\
&=&
       C_1\int_{\mathbb{A}}
          (V_{1}^{-1}(x))^{\frac{k_{1}r}{2}}
          (V_{1}(x))^{\frac{k_{1}r}{2}}
          |u_{m}|^{k_{1}r}dx
\nonumber\\
&&
      +C_1\int_{\mathbb{A}}
       \left(V_{1}(x)|u_{m}|^{2}\right)^{\frac{k_{1}r-1}{2}}
       \left(V_{1}(x)|u_n|^{2}\right)^{\frac{1}{2}}
          (V_{1}^{-1}(x))^{\frac{k_{1}r}{2}}
          dx
\nonumber\\
&&
      +C_1\int_{\mathbb{A}}
       \left(V_{2}(x)|v_m|^{2}\right)^{\frac{k_{2}(k_{1}r-1)}{2k_{1}}}
       \left(V_{1}(x)|u_{m}|^{2}\right)^{\frac{1}{2}}
       (V_{2}^{-1}(x))^{\frac{k_{2}(k_{1}r-1)}{2k_{1}}}
       (V_{1}^{-1}(x))^{\frac{1}{2}}
          dx
\nonumber\\
&&
      +C_1\int_{\mathbb{A}}
       \left(V_{2}(x)|v_m|^{2}\right)^{\frac{k_{2}(k_{1}r-1)}{2k_{1}}}
       \left(V_{1}(x)|u_n|^{2}\right)^{\frac{1}{2}}
       (V_{2}^{-1}(x))^{\frac{k_{2}(k_{1}r-1)}{2k_{1}}}
       (V_{1}^{-1}(x))^{\frac{1}{2}}
          dx
\nonumber\\
&\leq&
        C_1
        \left(
        \int_{\mathbb{A}}
        (V_{1}^{-1}(x))^{\frac{k_{1}r}{2}\cdot\frac{2}{2-k_{1}r}}
         dx
        \right)^{\frac{2-k_{1}r}{2}}
        \left(
        \int_{\mathbb{A}}
        (V_{1}(x))^{\frac{k_{1}r}{2}\cdot\frac{2}{k_{1}r}}
        |u_{m}|^{k_{1}r\cdot\frac{2}{k_{1}r}}
        dx
        \right)^{\frac{k_{1}r}{2}}
\nonumber\\
&&
       +C_1
       \left(
       \int_{\mathbb{A}}
       (V_{1}(x)|u_{m}|^{2})^{\frac{k_{1}r-1}{2}\cdot\frac{2}{ k_{1}r-1}}
        dx
       \right)
       ^{\frac{k_{1}r-1}{2}}
       \left(
       \int_{\mathbb{A}}
       (V_{1}(x)|u_n|^{2})^{\frac{1}{2}\cdot2}
       dx
       \right)^{\frac{1}{2}}
        \left(
       \int_{\mathbb{A}}
      (V_{1}^{-1}(x))^{\frac{k_{1}r}{2}\cdot\frac{2}{2-k_{1}r} }
       dx
       \right)^{\frac{2-k_{1}r}{2}}
\nonumber\\
&&
       +C_1
       \left(
       \int_{\mathbb{A}}
       (V_{2}(x)|v_m|^{2})^{\frac{k_{2}(k_{1}r-1)}{2k_{1}}\cdot\frac{2k_{1}}{k_{2}(k_{1}r-1)}}
        dx
       \right)
       ^{\frac{k_{2}(k_{1}r-1)}{2k_{1}}}
       \left(
       \int_{\mathbb{A}}
       (V_{1}(x)|u_{m}|^{2})^{\frac{1}{2}\cdot2}
       dx
       \right)^{\frac{1}{2}}
\nonumber\\
&&
       \cdot \left(
       \int_{\mathbb{A}}
       (V_{2}^{-1}(x))^{\frac{k_{2}(k_{1}r-1)}{2k_{1}}\cdot\frac{2k_{1}}{k_{1}-k_{2}(k_{1}r-1)}}
       (V_{1}^{-1}(x))^{\frac{1}{2}\cdot\frac{2k_{1}}{k_{1}-k_{2}(k_{1}r-1)}}
       dx
       \right)^{\frac{k_{1}-k_{2}(k_{1}r-1)}{2k_{1}}}
\nonumber\\
&&
       +C_1
       \left(
       \int_{\mathbb{A}}
       (V_{2}(x)|v_m|^{2})^{\frac{k_{2}(k_{1}r-1)}{2k_{1}}\cdot\frac{2k_{1}}{k_{2}(k_{1}r-1)}}
        dx
       \right)
       ^{\frac{k_{2}(k_{1}r-1)}{2k_{1}}}
       \left(
       \int_{\mathbb{A}}
       (V_{1}(x)|u_n|^{2})^{\frac{1}{2}\cdot2}
       dx
       \right)^{\frac{1}{2}}
\nonumber\\
&&
       \cdot \left(
       \int_{\mathbb{A}}
       (V_{2}^{-1}(x))^{\frac{k_{2}(k_{1}r-1)}{2k_{1}}\cdot\frac{2k_{1}}{k_{1}-k_{2}(k_{1}r-1)}}
       (V_{1}^{-1}(x))^{\frac{1}{2}\cdot\frac{2k_{1}}{k_{1}-k_{2}(k_{1}r-1)}}
       dx
       \right)^{\frac{k_{1}-k_{2}(k_{1}r-1)}{2k_{1}}}
\nonumber\\
&=&
       C_1
        \left(
        \int_{\mathbb{A}}
        (V_{1}^{-1}(x))^{\frac{k_{1}r}{2-k_{1}r}}
        dx
        \right)^{\frac{2-k_{1}r}{k_{1}r}\cdot\frac{k_{1}r}{2}}
        \left(
        \int_{\mathbb{A}}V_{1}(x)|u_{m}|^{2}dx
        \right)^{\frac{k_{1}r}{2}}
\nonumber\\
&&
       +C_1
       \left(\int_{\mathbb{A}}V_{1}(x)|u_{m}|^{2}dx\right)^{\frac{k_{1}r-1}{2}}
       \left(\int_{\mathbb{A}}V_{1}(x)|u_n|^{2}dx\right)^{\frac{1}{2}}
        \left(
        \int_{\mathbb{A}}
      (V_{1}^{-1}(x))^{\frac{k_{1}r}{2-k_{1}r}}
       dx
       \right)^{\frac{2-k_{1}r}{k_{1}r}\cdot\frac{k_{1}r}{2}}
\nonumber\\
&&
       +C_1
       \left(\int_{\mathbb{A}}V_{2}(x)|v_m|^{2}dx\right)^{\frac{k_{2}(k_{1}r-1)}{2k_{1}}}
       \left(\int_{\mathbb{A}}V_{1}(x)|u_{m}|^{2}dx\right)^{\frac{1}{2}}
       \left(
       \int_{\mathbb{A}}
       (V_{2}^{-1}(x))^{\frac{k_{2}(k_{1}r-1)}{k_{1}-k_{2}(k_{1}r-1)}}
       (V_{1}^{-1}(x))^{\frac{k_{1}}{k_{1}-k_{2}(k_{1}r-1)}}
       dx
       \right)^{\frac{k_{1}-k_{2}(k_{1}r-1)}{2k_{1}}}
\nonumber\\
&&
       +C_1
       \left(\int_{\mathbb{A}}V_{2}(x)|v_m|^{2}dx\right)^{\frac{k_{2}(k_{1}r-1)}{2k_{1}}}
       \left(\int_{\mathbb{A}}V_{1}(x)|u_n|^{2}dx\right)^{\frac{1}{2}}
       \left(
       \int_{\mathbb{A}}
       (V_{2}^{-1}(x))^{\frac{k_{2}(k_{1}r-1)}{k_{1}-k_{2}(k_{1}r-1)}}
       (V_{1}^{-1}(x))^{\frac{k_{1}}{k_{1}-k_{2}(k_{1}r-1)}}
       dx
       \right)^{\frac{k_{1}-k_{2}(k_{1}r-1)}{2k_{1}}}
\nonumber\\
&\leq&
       C_1
       \left(
        \int_{\mathbb{A}}
        (V_{1}^{-1}(x))^{\frac{k_{1}r}{2-k_{1}r}}
        dx
        \right)^{\frac{2-k_{1}r}{k_{1}r}\cdot\frac{k_{1}r}{2}}
        \left(
        \int_{\mathbb{A}}(|\nabla u_{m}|^{2}+V_{1}(x)|u_{m}|^{2})dx
        \right)^{\frac{k_{1}r}{2}}
\nonumber\\
&&
       +C_1
       \left(\int_{\mathbb{A}}(|\nabla u_{m}|^{2}+V_{1}(x)|u_{m}|^{2})dx\right)^{\frac{k_{1}r-1}{2}}
       \left(\int_{\mathbb{A}}(|\nabla u_n|^{2}+V_{1}(x)|u_n|^{2})dx\right)^{\frac{1}{2}}
       \left(
        \int_{\mathbb{A}}
      (V_{1}^{-1}(x))^{\frac{k_{1}r}{2-k_{1}r}}
       dx
       \right)^{\frac{2-k_{1}r}{k_{1}r}\cdot\frac{k_{1}r}{2}}
\nonumber\\
&&
       +C_1
       \left(\int_{\mathbb{A}}(|\nabla v_m|^{2}+V_{2}(x)|v_m|^{2})dx\right)
       ^{\frac{k_{2}(k_{1}r-1)}{2k_{1}}}
       \left(\int_{\mathbb{A}}(|\nabla u_{m}|^{2}+V_{1}(x)|u_{m}|^{2})dx\right)
       ^{\frac{1}{2}}
\nonumber\\
&&
       \cdot \left(
       \int_{\mathbb{A}}
       (V_{2}^{-1}(x))^{\frac{k_{2}(k_{1}r-1)}{k_{1}-k_{2}(k_{1}r-1)}}
       (V_{1}^{-1}(x))^{\frac{k_{1}}{k_{1}-k_{2}(k_{1}r-1)}}
       dx
       \right)^{\frac{k_{1}-k_{2}(k_{1}r-1)}{2k_{1}}}
\nonumber\\
&&
       +C_1
       \left(\int_{\mathbb{A}}(|\nabla v_m|^{2}+V_{2}(x)|v_m|^{2})dx\right)
       ^{\frac{k_{2}(k_{1}r-1)}{2k_{1}}}
       \left(\int_{\mathbb{A}}(|\nabla u_n|^{2}+V_{1}(x)|u_n|^{2})dx\right)
       ^{\frac{1}{2}}
\nonumber\\
&&
       \cdot \left(
       \int_{\mathbb{A}}
       (V_{2}^{-1}(x))^{\frac{k_{2}(k_{1}r-1)}{k_{1}-k_{2}(k_{1}r-1)}}
       (V_{1}^{-1}(x))^{\frac{k_{1}}{k_{1}-k_{2}(k_{1}r-1)}}
       dx
       \right)^{\frac{k_{1}-k_{2}(k_{1}r-1)}{2k_{1}}}
\nonumber\\
&=&
        C_1
        \left(
        \int_{\mathbb{A}}
        (V_{1}^{-1}(x))^{\frac{k_{1}r}{2-k_{1}r}}
        dx
        \right)^{\frac{2-k_{1}r}{k_{1}r}\cdot\frac{k_{1}r}{2}}
        \|u_{m}\|_{W_{1}}^{k_{1}r}
\nonumber\\
&&
       +C_1
       \left(
        \int_{\mathbb{A}}
        (V_{1}^{-1}(x))^{\frac{k_{1}r}{2-k_{1}r}}
         dx
       \right)^{\frac{2-k_{1}r}{k_{1}r}\cdot\frac{k_{1}r}{2}}
       \|u_{m}\|_{W_{1}}^{k_{1}r-1}
       \|u_n\|_{W_{1}}
\nonumber\\
&&
       +C_1
       \left(
       \int_{\mathbb{A}}
       (V_{2}^{-1}(x))^{\frac{k_{2}(k_{1}r-1)}{k_{1}-k_{2}(k_{1}r-1)}}
       (V_{1}^{-1}(x))^{\frac{k_{1}}{k_{1}-k_{2}(k_{1}r-1)}}
       dx
       \right)^{\frac{k_{1}-k_{2}(k_{1}r-1)}{2k_{1}}}
        \|v_m\|_{W_{2}}^{\frac{k_{2}(k_{1}r-1)}{k_{1}}}
        \|u_{m}\|_{W_{1}}
\nonumber\\
&&
       +C_1
       \left(
       \int_{\mathbb{A}}
       (V_{2}^{-1}(x))^{\frac{k_{2}(k_{1}r-1)}{k_{1}-k_{2}(k_{1}r-1)}}
       (V_{1}^{-1}(x))^{\frac{k_{1}}{k_{1}-k_{2}(k_{1}r-1)}}
       dx
       \right)^{\frac{k_{1}-k_{2}(k_{1}r-1)}{2k_{1}}}
        \|v_m\|_{W_{2}}^{\frac{k_{2}(k_{1}r-1)}{k_{1}}}
        \|u_n\|_{W_{1}}
\nonumber\\
&<&
       \epsilon^{\frac{k_{1}r}{2}} C_1 \|u_{m}\|_{W_{1}}^{k_{1}r}
       +\epsilon^{\frac{k_{1}r}{2}} C_1\|u_{m}\|_{W_{1}}^{k_{1}r-1}\|u_n\|_{W_{1}}
\nonumber\\
&&
       +C_1 C
       \left(
          \epsilon^{\frac{k_{2}r(k_{1}-k_{2}(k_{1}r-1))}{2k_{1}(2-k_{2}r)}}
          +\epsilon^{\frac{k_{1}r(k_{1}-k_{2}(k_{1}r-1))}{2k_{1}(2-k_{1}r)}}
       \right)
       \left(
       \|v_m\|_{W_{2}}^{\frac{k_{2}(k_{1}r-1)}{k_{1}}}\|u_{m}\|_{W_{1}}
      +\|v_m\|_{W_{2}}^{\frac{k_{2}(k_{1}r-1)}{k_{1}}}\|u_n\|_{W_{1}}
       \right).
\end{eqnarray*}
As $\{u_{n}\}$ and $\{v_m\}$ are bounded in $W_{1}$ and $W_{2}$, respectively,  the estimate just above indicates that there exist $M_{1}>0$ such that
\begin{eqnarray}\label{3.2.4.1}
\int_{\mathbb{A}}|\bar{F}_{t}(x,u_{m},v_m)||u_{m}-u_n|dx
<M_{1}\left(  \epsilon^{\frac{k_{1}r}{2}}
          +\epsilon^{\frac{k_{2}r(k_{1}-k_{2}(k_{1}r-1))}{2k_{1}(2-k_{2}r)}}
          +\epsilon^{\frac{k_{1}r(k_{1}-k_{2}(k_{1}r-1))}{2k_{1}(2-k_{1}r)}}
       \right)
\;\;\mbox{for any $n,m \in \mathbb{N}$}.
\end{eqnarray}
Moreover, by $u_{n}\rightharpoonup u$ in $W_{1}$, it is easy to see that  $u_n|_{B_{R}}\rightarrow u|_{B_{R}}$ in $L^{k_1}(B_{R})$ from Remark \ref{remark2.1} and $k_1\geq 2$. So that $\{u_{n}\}$ is a Cauchy sequence in $L^{k_1}(B_{R})$, i.e., for $\epsilon>0$ above, there exists $n_{\epsilon}>0$ such that
\begin{eqnarray*}
\left(\int_{B_{R}(0)}|u_m-u_n|^{k_{1}}dx\right)^{\frac{1}{k_{1}}}<\epsilon,\;\;\forall\; n,m>n_{\epsilon}.
\end{eqnarray*}
Then for all $n,m>n_{\epsilon}$, by $(F'_{1})$, $2\leq k_{i}<2^{\ast}$($i=1,2$), (\ref{A3}),
H\"{o}lder's inequality, Young's inequality and Remark \ref{remark2.1}, we have
\begin{eqnarray*}
&&
       \int_{B_{R}(0)}\left|\bar{F}_{t}(x,u_m,v_m)\right||u_m-u_n|dx
\nonumber\\
&&
\leq
        C_{1}\int_{B_{R}(0)}
        \left(
        |u_m|^{k_{1}r-1}+ |v_m|^{\frac{k_{2}(k_{1}r-1)}{k_{1}}}
        \right)
        |u_m-u_n|dx
\nonumber\\
&&
\leq
        C_{1}
        \left(
        \int_{B_{R}(0)}
        \left(
        |u_m|^{k_{1}r-1}+|v_m|^{\frac{k_{2}(k_{1}r-1)}{k_{1}}}
        \right)^{\frac{k_{1}}{k_{1}-1}}
        dx
        \right)^{\frac{k_{1}-1}{k_{1}}}
        \left(\int_{B_{R}(0)}|u_m-u_n|^{k_{1}}dx\right)^{\frac{1}{k_{1}}}
        \nonumber\\
 &&
 \leq
       \epsilon C_{1}
       \left(
       \int_{B_{R}(0)}
        \left(
        |u_m|^{k_{1}r-1}+|v_m|^{\frac{k_{2}(k_{1}r-1)}{k_{1}}}
        \right)^{\frac{k_{1}}{k_{1}-1}}
        dx
        \right)^{\frac{k_{1}-1}{k_{1}}}
        \nonumber\\
&&
 \leq   \epsilon C_{1}2^{\frac{1}{k_{1}}}
        \left(
        \int_{B_{R}(0)}
       \left(
       |u_m|^{\frac{k_{1}(k_{1}r-1)}{k_{1}-1}}+|v_m|^{\frac{k_{2}(k_{1}r-1)}{k_{1}-1}}
       \right)
       dx
       \right)^{\frac{k_{1}-1}{k_{1}}}
        \nonumber\\
&&
 \leq   \epsilon C_{1}2^{\frac{1}{k_{1}}}
        \left(
        \int_{B_{R}(0)}
       \left(
       \frac{k_{1}r-1}{k_{1}-1}|u_m|^{\frac{k_{1}(k_{1}r-1)}{k_{1}-1}\cdot\frac{k_{1}-1}{k_{1}r-1}}
       +\frac{k_{1}-k_{1}r}{ k_{1}-1}
       +\frac{k_{1}r-1}{k_{1}-1}|v_m|^{\frac{k_{2}(k_{1}r-1)}{k_{1}-1}\cdot\frac{k_{1}-1}{k_{1}r-1}}
       +\frac{k_{1}-k_{1}r}{ k_{1}-1}
       \right)dx
       \right)^{\frac{k_{1}-1}{k_{1}}}
        \nonumber\\
&&
 \leq
        \epsilon C_{1}2^{\frac{1}{k_{1}}}
        \left(
        \int_{B_{R}(0)}
        \left(
        |u_m|^{k_{1}}+ |v_m|^{k_{2}}+2
        \right)
        dx
        \right)^{\frac{k_{1}-1}{k_{1}}}
\nonumber\\
&&
=
        \epsilon C_{1} 2^{\frac{1}{k_{1}}}
            \left(
                   \|u_m\|_{k_{1}}^{k_{1}}
                 + \|v_m\|_{k_{2}}^{k_{2}}
                 + 2\mu(B_{R}(0))
           \right)^{\frac{k_{1}-1}{k_{1}}},
\end{eqnarray*}
where $\mu(B_{R}(0))$ denotes the Lebesgue measure of $B_{R}(0)$.
Combining the boundedness of $\|u_m\|_{k_{1}}$, $\|v_m\|_{k_{2}}$ and $\mu(B_{R}(0))$,
there exist $M_{2}>0$ such that
\begin{eqnarray}\label{3.2.4.2}
\int_{B_{R}(0)}|\bar{F}_{t}(x,u_{m},v_m)||u_{m}-u_n|dx
<\epsilon M_{2}
\;\;\mbox{for any $n,m>n_{\epsilon}$}.
\end{eqnarray}
Now, by \eqref{3.2.4.1} and \eqref{3.2.4.2}, we mention that
\begin{eqnarray*}
 \int_{\mathbb{R}^N}\left|\bar{F}_{t}(x,u_m,v_m)\right||u_m-u_n|dx
<M_{1}\left(  \epsilon^{\frac{k_{1}r}{2}}
          +\epsilon^{\frac{k_{2}r(k_{1}-k_{2}(k_{1}r-1))}{2k_{1}(2-k_{2}r)}}
          +\epsilon^{\frac{k_{1}r(k_{1}-k_{2}(k_{1}r-1))}{2k_{1}(2-k_{1}r)}}
     \right)
+\epsilon M_{2}
\;\;\mbox{for any $n,m>n_{\epsilon}$}.
\end{eqnarray*}
Since $\epsilon$ is arbitrary, it holds
\begin{eqnarray*}
 \int_{\mathbb{R}^N}\left|\bar{F}_{t}(x,u_m,v_m)\right||u_m-u_n|dx\rightarrow0 \;\;\mbox{as}\;\;n,m\rightarrow+\infty.
\end{eqnarray*}
Similarly, by $(\Lambda_{1})$, $(F_{1}')$, Remark \ref{remark2.1}, (\ref{A3})--(\ref{A5}) and $2\leq k_{i}<2^{\ast}$($i=1,2$),
we can obtain that
\begin{eqnarray*}
   \int_{\mathbb{R}^N}\left|\bar{F}_{s}(x,u_m,v_m)\right||v_m-v_n|dx\rightarrow0 \;\;\mbox{as}\;\;n,m\rightarrow+\infty.
\end{eqnarray*}
The proof is completed.
\qed

\par
\noindent
\begin{lemma}\label{lemma3.3}
Assume that $(F'_{5})$ and $( K)$ hold.
If there exists a sequence  $(u_{n},v_{n}) \rightharpoonup (u,v)$ in $W$ as $n\rightarrow +\infty$, then
\begin{eqnarray*}
    \int_{\mathbb{R}^N}k(x)\left|\bar{G}_{t}(u_m,v_m)\right||u_m-u_n|dx\rightarrow 0
\;\mbox{and}\;
    \int_{\mathbb{R}^N}k(x)\left|\bar{G}_{s}(u_m,v_m)\right||v_m-v_n|dx\rightarrow 0
\end{eqnarray*}
as $n,m\rightarrow+\infty$.
\end{lemma}
\par
\noindent
{\bf Proof.} The proof of this Lemma is similar to the proof of Lemma 3.5 in \cite{Liu-Zhang-Wang2024}.
\qed
\par
\noindent
\begin{lemma}\label{lemma3.5}
Assume that $(\Upsilon_{1})$, $(\Upsilon_{2})$, $(\Lambda_{0})$, $(\Lambda_{1})$, $(F'_{1})$, $(F'_{5})$ and $( K)$ hold.
For any given $c\in\R$ and $\varepsilon\in[0,1]$, $\bar{J}_{\varepsilon}$ satisfies $(PS)_{c}$-condition uniformly on $\varepsilon$.
\end{lemma}
\par
\noindent
{\bf Proof.}
Let  $c\in\R$ be fixed and $\{(\varepsilon_{n},(u_{n},v_{n}))\}\subset [0,1]\times W$ be any sequence such that
\begin{equation}\label{3.5.0}
     \bar{J}_{\varepsilon_{n}}(u_{n},v_{n})\rightarrow c
\;\mbox{and}\;
     \bar{J}_{\varepsilon_{n}}^{'}(u_{n},v_{n})\rightarrow 0
\;\mbox{as}\; n\rightarrow \infty.
\end{equation}
Then by \eqref{3.4.2}, $\{\varepsilon_{n}\}$ and $\{(u_{n},v_{n})\}$ are bounded.
Therefore, there exists a subsequence of $\{\varepsilon_{n}\}$, still denoted by $\{\varepsilon_{n}\}$, such that $\{\varepsilon_{n}\}$ converges to $\varepsilon$, and a subsequence of $\{(u_{n},v_{n})\}$, still denoted by $\{(u_{n},v_{n})\}$,  such that $\{(u_{n},v_{n})\}$ converges to $(u,v)$ weakly in $W$ and a.e. on $\R^{N}$.
Next, we shall show that this convergence is a strong convergence actually.
Evidently, it suffices to show that
\begin{eqnarray*}
u_{n}\rightarrow u \;\mbox{in} \;W_{1}
\;\mbox{and} \;
v_{n}\rightarrow v \;\mbox{in} \;W_{2}
\end{eqnarray*}
as $n\rightarrow\infty$.
Since $(W_{i},\|\cdot\|_{W_{i}})$($i=1,2$) are Banach spaces, it suffices to prove that $\{u_{n}\}$ is a Cauchy sequence in $W_{1}$ and $\{v_{n}\}$ is a Cauchy sequence in $W_{2}$.
Let
\begin{eqnarray*}
  \Gamma(u,v)
:= \int_{\mathbb{R}^{N}}\Phi_{1}\left(\frac{V_{1}(x)u^{2}+|\nabla u|^{2}}{2}\right)dx
   +\int_{\mathbb{R}^{N}}\Phi_{2}\left(\frac{V_{2}(x)v^{2}+|\nabla v|^{2}}{2}\right)dx,
\;\;\forall (u,v)\in W.
\end{eqnarray*}
It follows that
\begin{eqnarray}\label{3.5.2}
   \langle \Gamma'(u,v), (\varphi_{1},\varphi_{2})\rangle
&=&
\int_{\mathbb{R}^{N}}\phi_{1}\left(\frac{V_{1}(x)u^{2}+|\nabla u|^{2}}{2}\right)(V_{1}(x)u\varphi_{1}+\nabla u\cdot \nabla \varphi_{1} )dx
\nonumber\\
&&
+\int_{\mathbb{R}^{N}}\phi_{2}\left(\frac{V_{2}(x)v^{2}+|\nabla v|^{2}}{2}\right)(V_{2}(x)v\varphi_{2}+\nabla v\cdot \nabla \varphi_{2} )dx,
\;\;\forall (\varphi_{1},\varphi_{2})\in W.
\end{eqnarray}
For $(w_{1},h_{1}), (w_{2},h_{2}) \in W$,
let
$e_{1}=((V_{1}(x))^{\frac{1}{2}}w_{1},\nabla w_{1})$,
$e_{2}=((V_{1}(x))^{\frac{1}{2}}w_{2},\nabla w_{2})$,
$e_{3}=((V_{2}(x))^{\frac{1}{2}}h_{1},\nabla h_{1})$
and
$e_{4}=((V_{2}(x))^{\frac{1}{2}}h_{2},\nabla h_{2})$,
so that
\begin{eqnarray*}
e_{1}\cdot e_{2}=V_{1}(x)w_{1}w_{2}+\nabla w_{1}\cdot\nabla w_{2}
\;\;\mbox{and}\;\;
e_{3}\cdot e_{4}=V_{2}(x)h_{1}h_{2}+\nabla h_{1}\cdot\nabla h_{2}.
\end{eqnarray*}
Then by $(\Upsilon_{2})$, we observe that
\begin{eqnarray*}\label{3.5.3}
&&
\Gamma(w_{1},h_{1})-\Gamma(w_{2},h_{2})
\nonumber\\
&= &
\int_{\R^{N}}\left(
    \hbar_{1}\left(\frac{|e_{1}|}{\sqrt{2}}\right)
    -\hbar_{1}\left(\frac{|e_{2}|}{\sqrt{2}}\right)
            \right)dx
+\int_{\R^{N}}\left(
      \hbar_{2}\left(\frac{|e_{3}|}{\sqrt{2}}\right)
      -\hbar_{2}\left(\frac{|e_{4}|}{\sqrt{2}}\right)
            \right)dx
\nonumber\\
&\geq &
\int_{\R^{N}}\left(
      \hbar_{1}'\left(\frac{|e_{2}|}{\sqrt{2}}\right)\frac{|e_{1}|-|e_{2}|}{\sqrt{2}}
    +\frac{l_{1}}{2}(|e_{1}|-|e_{2}|)^{2}
            \right)dx
+\int_{\R^{N}}\left(
      \hbar_{2}'\left(\frac{|e_{4}|}{\sqrt{2}}\right)\frac{|e_{3}|-|e_{4}|}{\sqrt{2}}
    +\frac{l_{2}}{2}(|e_{3}|-|e_{4}|)^{2}
            \right)dx
\nonumber\\
&\geq &
\int_{\R^{N}}\phi_{1}\left(\frac{|e_{2}|^{2}}{2}\right)|e_{2}|\left(|e_{1}|-|e_{2}|\right)dx
+\frac{\min\{l_{1},\rho_{0}\}}{2}\int_{\R^{N}}(|e_{1}|-|e_{2}|)^{2}dx
\nonumber\\
&&
+\int_{\R^{N}}\phi_{2}\left(\frac{|e_{4}|^{2}}{2}\right)|e_{4}|\left(|e_{3}|-|e_{4}|\right)dx
+\frac{\min\{l_{2},\rho_{0}\}}{2}\int_{\R^{N}}(|e_{3}|-|e_{4}|)^{2}dx
\nonumber\\
&=&
\int_{\R^{N}}\phi_{1}\left(\frac{|e_{2}|^{2}}{2}\right)|e_{2}|\left(|e_{1}|-|e_{2}|\right)dx
+\frac{\min\{l_{1},\rho_{0}\}}{2}\int_{\R^{N}}(|e_{1}-e_{2}|^{2}+2e_{1}\cdot e_{2}-2|e_{1}||e_{2}|)dx
\nonumber\\
&&
+\int_{\R^{N}}\phi_{2}\left(\frac{|e_{4}|^{2}}{2}\right)|e_{4}|\left(|e_{3}|-|e_{4}|\right)dx
+\frac{\min\{l_{2},\rho_{0}\}}{2}\int_{\R^{N}}(|e_{3}-e_{4}|^{2}+2e_{3}\cdot e_{4}-2|e_{3}||e_{4}|)dx.
\end{eqnarray*}
We also  have
\begin{eqnarray*}\label{3.5.4}
&&
\langle \Gamma'(w_{2},h_{2}), (w_{1}-w_{2},h_{1}-h_{2})\rangle
\nonumber\\
&= &
\int_{\R^{N}}\phi_{1}\left(\frac{|e_{2}|^{2}}{2}\right)
(V_{1}(x)w_{2}(w_{1}-w_{2})+\nabla w_{2}\cdot \nabla (w_{1}-w_{2}))dx
\nonumber\\
&&
+\int_{\R^{N}}\phi_{2}\left(\frac{|e_{4}|^{2}}{2}\right)
(V_{2}(x)h_{2}(h_{1}-h_{2})+\nabla h_{2}\cdot \nabla (h_{1}-h_{2}))dx
\nonumber\\
&= &
\int_{\R^{N}}\phi_{1}\left(\frac{|e_{2}|^{2}}{2}\right)e_{2}\cdot (e_{1}-e_{2})dx
+\int_{\R^{N}}\phi_{2}\left(\frac{|e_{4}|^{2}}{2}\right)e_{4}\cdot (e_{3}-e_{4})dx.
\end{eqnarray*}
Hence, by $(\Upsilon_{1})$, we obtain
\begin{eqnarray}\label{3.5.5}
& &
    \Gamma(w_{1},h_{1})-\Gamma(w_{2},h_{2})-\langle \Gamma'(w_{2},h_{2}), (w_{1}-w_{2},h_{1}-h_{2})\rangle
\nonumber\\
&\geq&
    \int_{\R^{N}}\phi_{1}\left(\frac{|e_{2}|^{2}}{2}\right)|e_{2}|\left(|e_{1}|-|e_{2}|\right)dx
    +\frac{\min\{l_{1},\rho_{0}\}}{2}\int_{\R^{N}}(|e_{1}-e_{2}|^{2}+2e_{1}\cdot e_{2}-2|e_{1}||e_{2}|)dx
\nonumber\\
& &
     -\int_{\R^{N}}\phi_{1}\left(\frac{|e_{2}|^{2}}{2}\right)e_{2}\cdot (e_{1}-e_{2})dx
     + \int_{\R^{N}}\phi_{2}\left(\frac{|e_{4}|^{2}}{2}\right)|e_{4}|\left(|e_{3}|-|e_{4}|\right)dx
\nonumber\\
&&
+\frac{\min\{l_{2},\rho_{0}\}}{2}\int_{\R^{N}}(|e_{3}-e_{4}|^{2}+2e_{3}\cdot e_{4}-2|e_{3}||e_{4}|)dx
-\int_{\R^{N}}\phi_{2}\left(\frac{|e_{4}|^{2}}{2}\right)e_{4}\cdot (e_{3}-e_{4})dx
\nonumber\\
&= &
   \int_{\R^{N}}\left(\phi_{1}\left(\frac{|e_{2}|^{2}}{2}\right)-\min\{l_{1},\rho_{0}\} \right)
       (|e_{1}||e_{2}|-e_{1}\cdot e_{2})dx
   +\frac{\min\{l_{1},\rho_{0}\}}{2}\int_{\R^{N}}|e_{1}-e_{2}|^{2}dx
\nonumber\\
&&
+\int_{\R^{N}}\left(\phi_{2}\left(\frac{|e_{4}|^{2}}{2}\right)-\min\{l_{2},\rho_{0}\} \right)
       (|e_{3}||e_{4}|-e_{3}\cdot e_{4})dx
   +\frac{\min\{l_{2},\rho_{0}\}}{2}\int_{\R^{N}}|e_{3}-e_{4}|^{2}dx
\nonumber\\
&\geq &
   \int_{\R^{N}}\left(\rho_{0}-\min\{l_{1},\rho_{0}\} \right)(|e_{1}||e_{2}|-e_{1}\cdot e_{2})dx
   +\frac{\min\{l_{1},\rho_{0}\}}{2}\int_{\R^{N}}|e_{1}-e_{2}|^{2}dx
\nonumber\\
&&
+\int_{\R^{N}}\left(\rho_{0}-\min\{l_{2},\rho_{0}\} \right)(|e_{3}||e_{4}|-e_{3}\cdot e_{4})dx
   +\frac{\min\{l_{2},\rho_{0}\}}{2}\int_{\R^{N}}|e_{3}-e_{4}|^{2}dx
\nonumber\\
&\geq &
   \frac{\min\{l_{1},\rho_{0}\}}{2}\int_{\R^{N}}|e_{1}-e_{2}|^{2}dx
   +\frac{\min\{l_{2},\rho_{0}\}}{2}\int_{\R^{N}}|e_{3}-e_{4}|^{2}dx
\nonumber\\
&=&
   \frac{\min\{l_{1},\rho_{0}\}}{2}\|w_{1}-w_{2}\|_{W_{1}}^{2}
   +\frac{\min\{l_{2},\rho_{0}\}}{2}\|h_{1}-h_{2}\|_{W_{2}}^{2}.
\end{eqnarray}
Moreover, it follows that
\begin{eqnarray*}\label{3.5.6}
&&
    \bar{J}_{\varepsilon_{n}}(w_{1},h_{1})
   -\bar{J}_{\varepsilon_{m}}(w_{2},h_{2})
   -\langle \bar{J}_{\varepsilon_{m}}'(w_{2},h_{2}),(w_{1}-w_{2},h_{1}-h_{2})\rangle
\nonumber\\
&=&
     \Gamma(w_{1},h_{1})
     -\Gamma(w_{2},h_{2})
     -\langle \Gamma'(w_{2},h_{2}),(w_{1}-w_{2},h_{1}-h_{2})\rangle
\nonumber\\
&&
     +\int_{\mathbb{R}^{N}} \bar{F}_{t}(x,w_{2},h_{2})(w_{1}-w_{2})dx
     +\int_{\mathbb{R}^{N}}\bar{F}_{s}(x,w_{2},h_{2})(h_{1}-h_{2})dx
\nonumber\\
&&
   +\varepsilon_{m}\int_{\mathbb{R}^{N}}k(x) \bar{G}_{t}(w_{2},h_{2})(w_{1}-w_{2})dx
     +\varepsilon_{m}\int_{\mathbb{R}^{N}}k(x)
     \bar{G}_{s}(w_{2},h_{2})(h_{1}-h_{2})dx
\nonumber\\
&&
     +\int_{\mathbb{R}^{N}}(\bar{F}(x,w_{2},h_{2})-\bar{F}(x,w_{1},h_{1}))dx
     +\varepsilon_{m}\int_{\mathbb{R}^{N}}k(x)\bar{G}(w_{2},h_{2})dx
     -\varepsilon_{n}\int_{\mathbb{R}^{N}}k(x)\bar{G}(w_{1},h_{1}))dx
\nonumber\\
&=&
     \Gamma(w_{1},h_{1})
     -\Gamma(w_{2},h_{2})
     -\langle \Gamma'(w_{2},h_{2}),(w_{1}-w_{2},h_{1}-h_{2})\rangle
\nonumber\\
&&
     +\int_{\mathbb{R}^{N}} \bar{F}_{t}(x,w_{2},h_{2})(w_{1}-w_{2})dx
     +\int_{\mathbb{R}^{N}} \bar{F}_{s}(x,w_{2},h_{2})(h_{1}-h_{2})dx
\nonumber\\
&&
   +\varepsilon_{m}\int_{\mathbb{R}^{N}}k(x) \bar{G}_{t}(w_{2},h_{2})(w_{1}-w_{2})dx
   +\varepsilon_{m}\int_{\mathbb{R}^{N}}k(x) \bar{G}_{s}(w_{2},h_{2})(h_{1}-h_{2})dx
\nonumber\\
&&
     +\int_{\mathbb{R}^{N}}(\bar{F}(x,w_{2},h_{2})-\bar{F}(x,w_{1},h_{2}))dx
     +\int_{\mathbb{R}^{N}}(\bar{F}(x,w_{1},h_{2})-\bar{F}(x,w_{1},h_{1}))dx
\nonumber\\
&&
     +\varepsilon_{m}\int_{\mathbb{R}^{N}}k(x)(\bar{G}(w_{1},h_{2})-\bar{G}(w_{1},h_{1}))dx
     +\varepsilon_{m}\int_{\mathbb{R}^{N}}k(x)(\bar{G}(w_{2},h_{2})-\bar{G}(w_{1},h_{2}))dx
\nonumber\\
&&
     +(\varepsilon_{m}-\varepsilon_{n})\int_{\mathbb{R}^{N}}k(x)\bar{G}(w_{1},h_{1})dx.
\end{eqnarray*}
So, combining the equality above and (\ref{3.5.5}), and applying the mean value theorem, there exist $\theta_{i}\in \mathbb{R}$ with $0<\theta_{i}<1$($i=1,2,3,4$) such that
\begin{eqnarray*}\label{3.5.7}
& &
    \frac{\min\{l_{1},\rho_{0}\}}{2}\|w_{1}-w_{2}\|_{W_{1}}^{2}
    +\frac{\min\{l_{2},\rho_{0}\}}{2}\|h_{1}-h_{2}\|_{W_{2}}^{2}
\nonumber\\
&\leq&
    \bar{J}_{\varepsilon_{n}}(w_{1},h_{1})
   -\bar{J}_{\varepsilon_{m}}(w_{2},h_{2})
   +\|\bar{J}_{\varepsilon_{m}}'(w_{2},h_{2})\|_{W^{\ast}}\|(w_{1}-w_{2},h_{1}-h_{2})\|
\nonumber\\
&&
      +\int_{\mathbb{R}^{N}}| \bar{F}_{t}(x,w_{2},h_{2})||w_{1}-w_{2}|dx
     +\int_{\mathbb{R}^{N}} |\bar{F}_{s}(x,w_{2},h_{2})||h_{1}-h_{2}|dx
\nonumber\\
&&
   +|\varepsilon_{m}|\int_{\mathbb{R}^{N}}|k(x) \bar{G}_{t}(w_{2},h_{2})||w_{1}-w_{2}|dx
   +|\varepsilon_{m}|\int_{\mathbb{R}^{N}}|k(x) \bar{G}_{s}(w_{2},h_{2})||h_{1}-h_{2}|dx
\nonumber\\
&&
     +\int_{\mathbb{R}^{N}}|\bar{F}(x,w_{2},h_{2})-\bar{F}(x,w_{1},h_{2})|dx
     +\int_{\mathbb{R}^{N}}|\bar{F}(x,w_{1},h_{2})-\bar{F}(x,w_{1},h_{1})|dx
\nonumber\\
&&
     +|\varepsilon_{m}|\int_{\mathbb{R}^{N}}|k(x)||\bar{G}(w_{1},h_{2})-\bar{G}(w_{1},h_{1})|dx
     +|\varepsilon_{m}|\int_{\mathbb{R}^{N}}|k(x)||\bar{G}(w_{2},h_{2})-\bar{G}(w_{1},h_{2})|dx
\nonumber\\
&&
     +|\varepsilon_{m}-\varepsilon_{n}|\int_{\mathbb{R}^{N}}|k(x)||\bar{G}(w_{1},h_{1})|dx
\nonumber\\
&=&
    \bar{J}_{\varepsilon_{n}}(w_{1},h_{1})
   -\bar{J}_{\varepsilon_{n}}(w_{2},h_{2})
   +\|\bar{J}_{\varepsilon_{n}}'(w_{2},h_{2})\|_{W^{\ast}}\|(w_{1}-w_{2},h_{1}-h_{2})\|
\nonumber\\
&&
      +\int_{\mathbb{R}^{N}}| \bar{F}_{t}(x,w_{2},h_{2})||w_{1}-w_{2}|dx
     +\int_{\mathbb{R}^{N}} |\bar{F}_{s}(x,w_{2},h_{2})||h_{1}-h_{2}|dx
\nonumber\\
&&
   +|\varepsilon_{m}|\int_{\mathbb{R}^{N}}|k(x)
   \bar{G}_{t}(w_{2},h_{2})||w_{1}-w_{2}|dx
   +|\varepsilon_{m}|\int_{\mathbb{R}^{N}}|k(x) \bar{G}_{s}(w_{2},h_{2})||h_{1}-h_{2}|dx
\nonumber\\
&&
    +\int_{\mathbb{R}^{N}}| \bar{F}_{s}(x,w_{1},h_{1}+\theta_{1}(h_{2}-h_{1}))||h_{2}-h_{1}|dx
     +\int_{\mathbb{R}^{N}}| \bar{F}_{t}(x,w_{1}+\theta_{2}(w_{2}-w_{1}),h_{2})||w_{2}-w_{1}|dx
\nonumber\\
&&
   +|\varepsilon_{m}|\int_{\mathbb{R}^{N}}|k(x) \bar{G}_{t}(w_{1}+\theta_{3}(w_{2}-w_{1}),h_{2})||w_{2}-w_{1}|dx
   +|\varepsilon_{m}|\int_{\mathbb{R}^{N}}|k(x) \bar{G}_{s}(w_{1},h_{1}+\theta_{4}(h_{2}-h_{1}))||h_{2}-h_{1}|dx
\nonumber\\
&&
     +|\varepsilon_{m}-\varepsilon_{n}|\int_{\mathbb{R}^{N}}|k(x)||\bar{G}(w_{1},h_{1})|dx,
\end{eqnarray*}
where $W^{\ast}$ denotes the dual space of $W$.
Taking $w_{1}=u_{n}$, $w_{2}=u_{m}$, $h_{1}=v_{n}$ and $h_{2}=v_{m}$, we have
\begin{eqnarray}\label{3.5.8}
& &
    \frac{\min\{l_{1},\rho_{0}\}}{2}\|u_{n}-u_{m}\|_{W_{1}}^{2}
    + \frac{\min\{l_{2},\rho_{0}\}}{2}\|v_{n}-v_{m}\|_{W_{2}}^{2}
\nonumber\\
&\leq&
    \bar{J}_{\varepsilon_{n}}(u_{n},v_{n})
   -\bar{J}_{\varepsilon_{m}}(u_{m},v_{m})
   +\|\bar{J}_{\varepsilon_{m}}'(u_{m},v_{m})\|_{W^{\ast}}\|(u_{n}-u_{m},v_{n}-v_{m})\|
\nonumber\\
&&
      +\int_{\mathbb{R}^{N}}| \bar{F}_{t}(x,u_{m},v_{m})||u_{n}-u_{m}|dx
     +\int_{\mathbb{R}^{N}} |\bar{F}_{s}(x,u_{m},v_{m})||v_{n}-v_{m}|dx
\nonumber\\
&&
   +|\varepsilon_{m}|\int_{\mathbb{R}^{N}}|k(x)
   \bar{G}_{t}(u_{m},v_{m})||u_{n}-u_{m}|dx
\nonumber\\
&&
+|\varepsilon_{m}|\int_{\mathbb{R}^{N}}|k(x)\bar{G}_{s}(u_{m},v_{m})||v_{n}-v_{m}|dx
\nonumber\\
&&
    +\int_{\mathbb{R}^{N}}| \bar{F}_{s}(x,u_{n},v_{n}+\theta_{1}(v_{m}-v_{n}))||v_{m}-v_{n}|dx
\nonumber\\
&&
     +\int_{\mathbb{R}^{N}}| \bar{F}_{t}(x,u_{n}+\theta_{2}(u_{m}-u_{n}),v_{m})||u_{m}-u_{n}|dx
\nonumber\\
&&
   +|\varepsilon_{m}|\int_{\mathbb{R}^{N}}|k(x) \bar{G}_{t}(u_{n}+\theta_{3}(u_{m}-u_{n}),v_{m})||u_{m}-u_{n}|dx
\nonumber\\
&&
   +|\varepsilon_{m}|\int_{\mathbb{R}^{N}}|k(x) \bar{G}_{s}(u_{n},v_{n}+\theta_{4}(v_{m}-v_{n}))||v_{m}-v_{n}|dx
\nonumber\\
&&
     +|\varepsilon_{m}-\varepsilon_{n}|\int_{\mathbb{R}^{N}}|k(x)||\bar{G}(u_{n},v_{n})|dx.
\end{eqnarray}
From (\ref{3.5.0}) we have
\begin{eqnarray*}\label{3.5.9}
\lim_{m\rightarrow\infty}\|\bar{J}_{\varepsilon_{m}}'(u_{m},v_{m})\|_{W^{\ast}}=0.
\end{eqnarray*}
The weak convergence of $\{(u_n,v_n)\}$ in $W$ implies that $\{\|(u_n,v_n)\|\}$ is a bounded sequence.
Hence
\begin{eqnarray}\label{3.5.10}
\|\bar{J}_{\varepsilon_{m}}'(u_{m},v_{m})\|_{W^{\ast}}\|(u_{n}-u_{m},v_{n}-v_{m})\|\rightarrow0 \;\;\mbox{as}\;\;n,m\rightarrow+\infty.
\end{eqnarray}
By (\ref{3.5.0}), we also have
\begin{eqnarray}\label{3.5.11}
\bar{J}_{\varepsilon_{n}}(u_{n},v_{n})-\bar{J}_{\varepsilon_{m}}(u_{m},v_{m})\rightarrow0 \;\;\mbox{as}\;\;n,m\rightarrow+\infty.
\end{eqnarray}
By arguments as in  Lemma \ref{lemma3.2} and Lemma \ref{lemma3.3} and the fact that  $0<\theta_{i}<1$($i=1,2,3,4$), we can obtain that
\begin{eqnarray}
&&\label{3.5.12}
\int_{\mathbb{R}^{N}}| \bar{F}_{s}(x,u_{n},v_{n}+\theta_{1}(v_{m}-v_{n}))||v_{m}-v_{n}|dx
\rightarrow0 \;\;\mbox{as}\;\;n,m\rightarrow+\infty,
\\
&&\label{3.5.13}
\int_{\mathbb{R}^{N}}| \bar{F}_{t}(x,u_{n}+\theta_{2}(u_{m}-u_{n}),v_{m})||u_{m}-u_{n}|dx
\rightarrow0 \;\;\mbox{as}\;\;n,m\rightarrow+\infty,
\\
&&\label{3.5.14}
|\varepsilon_{m}|\int_{\mathbb{R}^{N}}|k(x) \bar{G}_{t}(u_{n}+\theta_{3}(u_{m}-u_{n}),v_{m})||u_{m}-u_{n}|dx
\rightarrow0 \;\;\mbox{as}\;\;n,m\rightarrow+\infty,
\\
&&\label{3.5.15}
|\varepsilon_{m}|\int_{\mathbb{R}^{N}}|k(x) \bar{G}_{s}(u_{n},v_{n}+\theta_{4}(v_{m}-v_{n}))||v_{m}-v_{n}|dx
\rightarrow0 \;\;\mbox{as}\;\;n,m\rightarrow+\infty.
\end{eqnarray}
In addition, $\varepsilon_{n}\rightarrow \varepsilon$ in $\R$ implies that $\{\varepsilon_{n}\}$ is a Cauchy sequence in $\R$, i.e., for $\epsilon>0$ given in (\ref{3.2.2}), there exists $n_{\epsilon}>0$ such that
\begin{eqnarray*}
|\varepsilon_{m}-\varepsilon_{n}|<\epsilon,\;\;\forall\; n,m>n_{\epsilon}.
\end{eqnarray*}
Let $\mathbb{B}=\{x\in \R^{N}: |(u_{n},v_{n})|\leq \delta\}$.
Then for all $n,m>n_{\epsilon}$, by $(F'_{5})$ and $(K)$, we have
\begin{eqnarray*}
|\varepsilon_{m}-\varepsilon_{n}|\int_{\mathbb{R}^{N}}|k(x)||\bar{G}(u_{n},v_{n})|dx
&=&
|\varepsilon_{m}-\varepsilon_{n}|\int_{\mathbb{B}}|k(x)||\bar{G}(u_{n},v_{n})|dx
\nonumber\\
&\leq&
|\varepsilon_{m}-\varepsilon_{n}|\max_{|(t,s)|\leq \delta}|\bar{G}(t,s)|\int_{\mathbb{B}}|k(x)|dx
\nonumber\\
&\leq&
|\varepsilon_{m}-\varepsilon_{n}|\max_{|(t,s)|\leq \delta}|\bar{G}(t,s)|\int_{\R^{N}}|k(x)|dx
\nonumber\\
&<&
\epsilon\max_{|(t,s)|\leq \delta}|\bar{G}(t,s)|\int_{\R^{N}}|k(x)|dx.
\end{eqnarray*}
Since $\epsilon$ is arbitrary, it holds
\begin{eqnarray}\label{3.5.16}
 |\varepsilon_{m}-\varepsilon_{n}|\int_{\mathbb{R}^{N}}|k(x)||\bar{G}(u_{n},v_{n})|dx\rightarrow0 \;\;\mbox{as}\;\;n,m\rightarrow+\infty.
\end{eqnarray}
As a consequence, by (\ref{3.5.8})--(\ref{3.5.16}), Lemma \ref{lemma3.2} and Lemma \ref{lemma3.3},
we can now conclude that $\{u_{n}\}$ is a Cauchy sequence in $W_{1}$ and $\{v_{n}\}$ is a Cauchy sequence in $W_{2}$.
The proof is completed.\qed
\par
\noindent
\begin{lemma}\label{lemma3.6}
Assume that  $(\Upsilon_1)$, $(\Upsilon_3)$, $(F_0')$, $(F_2')$--$(F'_{4})$, $(\Lambda_{0})$ and $(\Lambda_{1})$ hold.
For any given $(u,v)\in W \backslash \{(0,0)\}$, $\bar{J}_{0}(u,v)$  satisfies condition $(A_{5})$ in Lemma \ref{lemma2.2}.
\end{lemma}
\par
\noindent
{\bf Proof.}
Let $(u,v)\in W \backslash \{(0,0)\}$ be fixed.
Note that  $\bar{F}$ is even on $(u,v)\in \mathbb{R} \times \mathbb{R}$. Without loss of generality, we may assume $\vartheta>0$ and consider
 \begin{eqnarray}\label{3.6.0}
&&
         \gamma_{(u,v)}(\vartheta)
=
         \vartheta^{-\beta}\bar{J}_{0}(\vartheta u,\vartheta v)
=
       \vartheta^{-\beta}\int_{\mathbb{R}^{N}}
        \left\{
        \Phi_{1}\left(\frac{V_{1}(x)u^{2}+|\nabla u|^{2}}{2}\vartheta^{2}\right)
        +\Phi_{2}\left(\frac{V_{2}(x)v^{2}+|\nabla v|^{2}}{2}\vartheta^{2}\right)
        \right\}dx
\nonumber\\
&&\qquad\qquad\qquad\qquad\qquad\qquad
          -\vartheta^{-\beta}\int_{\mathbb{R}^{N}}
          \bar{F}(x,\vartheta u,\vartheta v)dx.
\end{eqnarray}
It follows from $(\Upsilon_3)$, $(\Upsilon_1)$, $(F_2')$, $(\Lambda_{0})$ and $\beta\in (0,2)$ that
\begin{eqnarray}\label{3.6.1}
&&
          \gamma_{(u,v)}'(\vartheta)
\nonumber\\
&=&
         \vartheta^{-\beta-1}
          \int_{\mathbb{R}^{N}}
          \left\{
          \vartheta^{2}
          \phi_{1}\left(\frac{V_{1}(x)u^{2}+|\nabla u|^{2}}{2}\vartheta^{2}\right)
                (V_{1}(x)u^{2}+|\nabla u|^{2})
          -\beta\Phi_{1}\left(\frac{V_{1}(x)u^{2}+|\nabla u|^{2}}{2}\vartheta^{2}\right)
          \right\}dx
\nonumber\\
& &
          +\vartheta^{-\beta-1}
          \int_{\mathbb{R}^{N}}
          \left\{
          \vartheta^{2}
          \phi_{2}\left(\frac{V_{2}(x)v^{2}+|\nabla v|^{2}}{2}\vartheta^{2}\right)
              (V_{2}(x)v^{2}+|\nabla v|^{2})
          -\beta\Phi_{2}\left(\frac{V_{2}(x)v^{2}+|\nabla v|^{2}}{2}\vartheta^{2}\right)
          \right\}dx
\nonumber\\
&&
          - \vartheta^{-\beta-1}
          \int_{\mathbb{R}^{N}}
          \left(
          \bar{F}_{t}(x,\vartheta u,\vartheta v)
                     \vartheta u
          + \bar{F}_{s}(x,\vartheta u,\vartheta v)
                    \vartheta v
          - \beta\bar{F}(x,\vartheta u,\vartheta v)
          \right)dx
\nonumber\\
&=&
          \vartheta^{-\beta-1}
          \int_{\mathbb{R}^{N}}
          \left\{
          2\phi_{1}\left(\frac{V_{1}(x)u^{2}+|\nabla u|^{2}}{2}\vartheta^{2}\right)
          \frac{V_{1}(x)u^{2}+|\nabla u|^{2}}{2}\vartheta^{2}
          -\beta\Phi_{1}\left(\frac{V_{1}(x)u^{2}+|\nabla u|^{2}}{2}\vartheta^{2}\right)
          \right\}dx
\nonumber\\
& &
          +\vartheta^{-\beta-1}
          \int_{\mathbb{R}^{N}}
          \left\{
          2\phi_{2}\left(\frac{V_{2}(x)v^{2}+|\nabla v|^{2}}{2}\vartheta^{2}\right)
          \frac{V_{2}(x)v^{2}+|\nabla v|^{2}}{2}\vartheta^{2}
          -\beta\Phi_{2}\left(\frac{V_{2}(x)v^{2}+|\nabla v|^{2}}{2}\vartheta^{2}\right)
          \right\}dx
\nonumber\\
&&
          - \vartheta^{-\beta-1}
          \int_{\mathbb{R}^{N}}
          \left(
        \bar{F}_{t}(x,\vartheta u,\vartheta v)\vartheta u
          +\bar{F}_{s}(x,\vartheta u,\vartheta v)\vartheta v
          - \beta\bar{F}(x,\vartheta u,\vartheta v)
          \right)dx
\nonumber\\
&\geq&
         \vartheta^{-\beta-1}
          (2-\beta)\int_{\mathbb{R}^{N}}
          \left\{
          \Phi_{1}\left(\frac{V_{1}(x)u^{2}+|\nabla u|^{2}}{2}\vartheta^{2}\right)
          +\Phi_{2}\left(\frac{V_{2}(x)v^{2}+|\nabla v|^{2}}{2}\vartheta^{2}\right)
          \right\}dx
\nonumber\\
&&
          -\vartheta^{-\beta-1}
          \int_{\mathbb{R}^{N}}
          \left(
         \bar{F}_{t}(x,\vartheta u,\vartheta v)\vartheta u
          + \bar{F}_{s}(x,\vartheta u,\vartheta v)\vartheta v
          - \beta\bar{F}(x,\vartheta u,\vartheta v)
          \right)dx
\nonumber\\
&\geq&
         \frac{\rho_{0}(2-\beta)}{2}\vartheta^{-\beta-1}
        \int_{\mathbb{R}^{N}}(V_{1}(x)u^{2}+|\nabla u|^{2})\vartheta^{2}dx
        +
        \frac{\rho_{0}(2-\beta)}{2}\vartheta^{-\beta-1}
        \int_{\mathbb{R}^{N}}(V_{2}(x)v^{2}+|\nabla v|^{2})\vartheta^{2}dx
\nonumber\\
&&
          -\vartheta^{-\beta-1}
          \int_{\mathbb{R}^{N}}
          \left(
          \bar{F}_{t}(x,\vartheta u,\vartheta v)\vartheta u
          +\bar{F}_{s}(x,\vartheta u,\vartheta v)\vartheta v
          - \beta\bar{F}(x,\vartheta u,\vartheta v)
          \right)dx
\nonumber\\
&\geq&
         \frac{\rho_{0}(2-\beta)}{2}\vartheta^{-\beta-1}
        \int_{\mathbb{R}^{N}}
        (V_{1}(x)(\vartheta u)^{2}+V_{2}(x)(\vartheta v)^{2})dx
\nonumber\\
&&
          -\vartheta^{-\beta-1}
          \int_{\mathbb{R}^{N}}
          \left(
          \bar{F}_{t}(x,\vartheta u,\vartheta v)\vartheta u
          +\bar{F}_{s}(x,\vartheta u,\vartheta v)\vartheta v
          - \beta\bar{F}(x,\vartheta u,\vartheta v)
          \right)dx
\nonumber\\
&=&
        \vartheta^{-\beta-1}
        \int_{\mathbb{R}^{N}}
        \left\{
        \frac{\rho_{0}(2-\beta)}{2}
        (V_{1}(x)(\vartheta u)^{2}+V_{2}(x)(\vartheta v)^{2})
        -
        \left(
          \bar{F}_{t}(x,\vartheta u,\vartheta v)\vartheta u
          +\bar{F}_{s}(x,\vartheta u,\vartheta v)\vartheta v
          - \beta\bar{F}(x,\vartheta u,\vartheta v)
          \right)
          \right\}
          dx
\nonumber\\
&> &   0.
\end{eqnarray}
\par
Next, take $\eta>0$ small enough such that
$$
\mu(K_{\eta})> 0\;\; \mbox{and}\;\; K_{\eta}:=\left\{x \in \mathbb{R}^{N}: \eta <|(u,v)|<\frac{1}{\eta}\right\},
$$
where $\mu$ denotes the Lebesgue measure of $\mathbb{R}^{N}$.
Due to $(\Upsilon_{1})$ and $\bar{F}(x,u,v)\geq 0$, we can estimate the function $\gamma_{(u,v)}(\vartheta)$ as follows:
\begin{eqnarray*}\label{3.6.2}
&&
         \gamma_{(u,v)}(\vartheta)
\nonumber\\
&=&
       \vartheta^{-\beta}\int_{\mathbb{R}^{N}}
        \left\{
        \Phi_{1}\left(\frac{V_{1}(x)u^{2}+|\nabla u|^{2}}{2}\vartheta^{2}\right)
        +\Phi_{2}\left(\frac{V_{2}(x)v^{2}+|\nabla v|^{2}}{2}\vartheta^{2}\right)
        \right\}dx
\nonumber\\
&&
          -\vartheta^{-\beta}\int_{\mathbb{R}^{N}}
                  \bar{F}(x,\vartheta u,\vartheta v)dx
\nonumber\\
&\leq&
         \frac{\rho_{1}}{2}\vartheta^{2-\beta}
         \int_{\mathbb{R}^{N}}
        (V_{1}(x)u^{2}+|\nabla u|^{2}
        +V_{2}(x)v^{2}+|\nabla v|^{2})dx
       - \vartheta^{-\beta}\int_{\mathbb{R}^{N}}
           \bar{F}(x,\vartheta u,\vartheta v)dx
\nonumber\\
&=&
        \frac{\rho_{1}}{2}\vartheta^{2-\beta}(\|u\|_{W_{1}}^{2}+\|v\|_{W_{2}}^{2})
       -\vartheta^{-\beta}\int_{\mathbb{R}^{N}}
                     \bar{F}(x,\vartheta u,\vartheta v)dx
\nonumber\\
&\leq&
        \frac{\rho_{1}}{2}\vartheta^{2-\beta}(\|u\|_{W_{1}}^{2}+\|v\|_{W_{2}}^{2})
       - \vartheta^{-\beta}
         \int_{K_{\eta}}
         \frac{\bar{F}(x,\vartheta u,\vartheta v)
                      (|u|^{2}+|v|^{2})^{\frac{\beta}{2}}}{|u|^{\beta}+|v|^{\beta}}dx
\nonumber\\
&=&
        \frac{\rho_{1}}{2}\vartheta^{2-\beta}(\|u\|_{W_{1}}^{2}+\|v\|_{W_{2}}^{2})
         -\int_{K_{\eta}}\frac{\bar{F}(x,\vartheta u,\vartheta v)
                        (|u|^{2}+|v|^{2})^{\frac{\beta}{2}}}
             {(|\vartheta u|^{\beta}+|\vartheta v|^{\beta})}dx
\nonumber\\
&\leq&
        \frac{\rho_{1}}{2}\vartheta^{2-\beta}(\|u\|_{W_{1}}^{2}+\|v\|_{W_{2}}^{2})
       -\eta^{\beta}
         \int_{K_{\eta}}\frac{\bar{F}(x,\vartheta u,\vartheta v)}
                     {(|\vartheta u|^{\beta}+|\vartheta v|^{\beta})}dx
\nonumber\\
&\leq&
        \frac{\rho_{1}}{2}\vartheta^{2-\beta}(\|u\|_{W_{1}}^{2}+\|v\|_{W_{2}}^{2})
       -\eta^{\beta}\mu(K_{\eta})
         \inf_{x\in K_{\eta}}
         \frac{\bar{F}(x,\vartheta u,\vartheta v)}
                  {(|\vartheta u|^{\beta}+|\vartheta v|^{\beta})}.
\end{eqnarray*}
So $\lim_{\vartheta\rightarrow0^{+}}\gamma_{(u,v)}(\vartheta)= -\infty$ since  $(F'_{3})$ and $\beta\in (0,2)$.
Hence,  $\gamma_{(u,v)}(\vartheta)<0$ for all $\vartheta>0$ small enough.
\par
In addition, by $(\Upsilon_{1})$, (\ref{r3.4}), $(\Lambda_{0})$, $(\Lambda_{1})$, $\beta\in (0,2)$, $1<k_{i}r<2$($i=1,2$) and H\"{o}lder's inequality, we have, for $t>0$ large enough,
\begin{eqnarray*}\label{3.6.3}
 &&
         \gamma_{(u,v)}(\vartheta)
 \nonumber\\
&=&
        \vartheta^{-\beta}\int_{\mathbb{R}^{N}}
        \left\{
        \Phi_{1}\left(\frac{V_{1}(x)u^{2}+|\nabla u|^{2}}{2}\vartheta^{2}\right)
        +\Phi_{2}\left(\frac{V_{2}(x)v^{2}+|\nabla v|^{2}}{2}\vartheta^{2}\right)
        \right\}dx
\nonumber\\
&&
          - \vartheta^{-\beta}
            \int_{\mathbb{R}^{N}}\bar{F}(x,\vartheta u,\vartheta v)dx
\nonumber\\
&\geq&
         \frac{\rho_{0}}{2}\vartheta^{2-\beta}
         \int_{\mathbb{R}^{N}}
        (V_{1}(x)u^{2}+|\nabla u|^{2}
        +V_{2}(x)v^{2}+|\nabla v|^{2})dx
       - \vartheta^{-\beta}
               \int_{\mathbb{R}^{N}}\bar{F}(x,\vartheta u,\vartheta v)dx
\nonumber\\
&=&
        \frac{\rho_{0}}{2}\vartheta^{2-\beta}(\|u\|_{W_{1}}^{2}+\|v\|_{W_{2}}^{2})
       -\vartheta^{-\beta}
                \int_{\mathbb{R}^{N}}\bar{F}(x,\vartheta u,\vartheta v)dx
\nonumber\\
&\geq&
        \frac{\rho_{0}}{2}\vartheta^{2-\beta}(\|u\|_{W_{1}}^{2}+\|v\|_{W_{2}}^{2})
       -\vartheta^{-\beta}
             \int_{\mathbb{R}^{N}}|\bar{F}(x,\vartheta u,\vartheta v)|dx
\nonumber\\
&\geq&
        \frac{\rho_{0}}{2}\vartheta^{2-\beta}(\|u\|_{W_{1}}^{2}+\|v\|_{W_{2}}^{2})
       - C_{3}\vartheta^{-\beta}
              \int_{\mathbb{R}^{N}}(|\vartheta u|^{k_{1}r}+|\vartheta v|^{k_{2}r})dx
\nonumber\\
&=&
        \frac{\rho_{0}}{2}\vartheta^{2-\beta}(\|u\|_{W_{1}}^{2}+\|v\|_{W_{2}}^{2})
       - C_{3}\vartheta^{k_{1}r-\beta}\int_{\mathbb{R}^{N}}|u|^{k_{1}r}dx
       - C_{3}\vartheta^{k_{2}r-\beta}\int_{\mathbb{R}^{N}}|v|^{k_{2}r}dx
\nonumber\\
&=&
       \frac{\rho_{0}}{2}\vartheta^{2-\beta}(\|u\|_{W_{1}}^{2}+\|v\|_{W_{2}}^{2})
       - C_{3}\vartheta^{k_{1}r-\beta}
       \int_{\mathbb{R}^{N}}
          (V_{1}^{-1}(x))^{\frac{k_{1}r}{2}}
          (V_{1}(x))^{\frac{k_{1}r}{2}}
          |u|^{k_{1}r}dx
\nonumber\\
&&
      -C_3\vartheta^{k_{2}r-\beta}
      \int_{\mathbb{R}^{N}}
          (V_{2}^{-1}(x))^{\frac{k_{2}r}{2}}
          (V_{2}(x))^{\frac{k_{2}r}{2}}
          |v|^{k_{2}r}dx
\nonumber\\
&\geq&
       \frac{\rho_{0}}{2}\vartheta^{2-\beta}(\|u\|_{W_{1}}^{2}+\|v\|_{W_{2}}^{2})
\nonumber\\
&&
       -C_{3}\vartheta^{k_{1}r-\beta}
        \left(
        \int_{\mathbb{R}^{N}}
        (V_{1}^{-1}(x))^{\frac{k_{1}r}{2}\cdot\frac{2}{2-k_{1}r}}
         dx
        \right)^{\frac{2-k_{1}r}{2}}
        \left(
        \int_{\mathbb{R}^{N}}
        (V_{1}(x))^{\frac{k_{1}r}{2}\cdot\frac{2}{k_{1}r}}
        |u|^{k_{1}r\cdot\frac{2}{k_{1}r}}
        dx
        \right)^{\frac{k_{1}r}{2}}
\nonumber\\
&&
       -C_3\vartheta^{k_{2}r-\beta}
        \left(
        \int_{\mathbb{R}^{N}}
        (V_{2}^{-1}(x))^{\frac{k_{2}r}{2}\cdot\frac{2}{2-k_{2}r}}
         dx
        \right)^{\frac{2-k_{2}r}{2}}
        \left(
        \int_{\mathbb{R}^{N}}
        (V_{2}(x))^{\frac{k_{2}r}{2}\cdot\frac{2}{k_{2}r}}
        |v|^{k_{2}r\cdot\frac{2}{k_{2}r}}
        dx
        \right)^{\frac{k_{2}r}{2}}
\nonumber\\
&=&
       \frac{\rho_{0}}{2}\vartheta^{2-\beta}(\|u\|_{W_{1}}^{2}+\|v\|_{W_{2}}^{2})
\nonumber\\
&&
       -C_3\vartheta^{k_{1}r-\beta}
        \left(
        \int_{\mathbb{R}^{N}}
        (V_{1}^{-1}(x))^{\frac{k_{1}r}{2-k_{1}r}}
        dx
        \right)^{\frac{2-k_{1}r}{k_{1}r}\cdot\frac{k_{1}r}{2}}
        \left(
        \int_{\mathbb{R}^{N}}V_{1}(x)|u|^{2}dx
        \right)^{\frac{k_{1}r}{2}}
\nonumber\\
&&
       -C_3\vartheta^{k_{2}r-\beta}
        \left(
        \int_{\mathbb{R}^{N}}
        (V_{2}^{-1}(x))^{\frac{k_{2}r}{2-k_{2}r}}
         dx
        \right)^{\frac{2-k_{2}r}{k_{2}r}\cdot\frac{k_{2}r}{2}}
        \left(
        \int_{\mathbb{R}^{N}}V_{2}(x)|v|^{2}dx
        \right)^{\frac{k_{2}r}{2}}
\nonumber\\
&\geq&
        \frac{\rho_{0}}{2}\vartheta^{2-\beta}(\|u\|_{W_{1}}^{2}+\|v\|_{W_{2}}^{2})
       -C_5\vartheta^{k_{1}r-\beta}
        \left(
        \int_{\mathbb{R}^{N}}(|\nabla u|^{2}+V_{1}(x)|u|^{2})dx
        \right)^{\frac{k_{1}r}{2}}
 \nonumber\\
&&
      -C_5\vartheta^{k_{2}r-\beta}
        \left(
        \int_{\mathbb{R}^{N}}(|\nabla v|^{2}+V_{2}(x)|v|^{2})dx
        \right)^{\frac{k_{2}r}{2}}
\nonumber\\
&=&
        \frac{\rho_{0}}{2}\vartheta^{2-\beta}(\|u\|_{W_{1}}^{2}+\|v\|_{W_{2}}^{2})
      -C_5\vartheta^{k_{1}r-\beta}\|u\|_{W_{1}}^{k_{1}r}
       -C_5\vartheta^{k_{2}r-\beta}\|v\|_{W_{2}}^{k_{2}r}
\nonumber\\
&>&  0,
\end{eqnarray*}
where $C_{5}$ is given in (\ref{3.4.3}).
\par
Consequently, combining with (\ref{3.6.1}), for fixed $(u,v)\in W \backslash \{(0,0)\}$, we get $\gamma_{(u,v)}(\vartheta)$ has a unique zero $\vartheta(u,v)$ such that
$\gamma_{(u,v)}(\vartheta)<0$, for $0<\vartheta<\vartheta(u,v)$ and  $\gamma_{(u,v)}(\vartheta)\geq 0$, for $\vartheta(u,v)\leq \vartheta$.
Furthermore, the above result holds for $\bar{J}_{0}(\vartheta u, \vartheta v)$, i.e., $\bar{J}_{0}(u,v)$  satisfies condition $(A_{5})$.
\qed
\par
Next, we shall give a crucial result, which will be used to estimate the energy of solutions
to ensure that the solutions of problem \eqref{mod1} is also the solutions of the original problem \eqref{eq1}.
\par
\noindent
\begin{lemma}\label{lemma3.8}
Assume that  $(\Upsilon_{1})$, $(\Lambda_{0})$, $(\Lambda_{1})$, $(F'_{1})$, $(F'_{5})$ and $(K)$ hold.  There exist positive constants $C^{*}$, $D^{*}$, $\zeta$ and $\vartheta$ such that  if $\bar{J}_{\varepsilon}'(u,v)=0$ with $|\varepsilon|\leq 1$, then
$\|u\|_{\infty}\leq C^{*} \|u\|_{2^{\ast}}^{\zeta}$ and $\|v\|_{\infty}\leq D^{*} \|v\|_{2^{\ast}}^{\vartheta}$.
\end{lemma}
\par
\noindent
{\bf Proof.}
Let $(u,v)\in W$ be a critical point of $\bar{J}_{\varepsilon}$.
In view of system \eqref{mod1}, for any $\psi=(\psi_{1},\psi_{2})\in  W$,
we obtain
\begin{eqnarray}\label{3.8.1}
         \int_{\mathbb{R}^{N}}
         \phi_{1}\left(\frac{V_{1}(x)u^{2}+|\nabla u|^{2}}{2}\right)
          (V_{1}(x)u \psi_{1}+\nabla u\cdot \nabla  \psi_{1})dx
=
        \int_{\mathbb{R}^{N}} \bar{F}_{t}(x,u,v) \psi_{1}dx
         +\varepsilon\int_{\mathbb{R}^{N}}k(x) \bar{G}_{t}(u,v) \psi_{1}dx
\end{eqnarray}
and
\begin{eqnarray}\label{3.8.2}
         \int_{\mathbb{R}^{N}}
         \phi_{2}\left(\frac{V_{2}(x)v^{2}+|\nabla v|^{2}}{2}\right)
          (V_{2}(x)v \psi_{2}+\nabla v\cdot \nabla  \psi_{2})dx
=
        \int_{\mathbb{R}^{N}} \bar{F}_{s}(x,u,v) \psi_{2}dx
         +\varepsilon\int_{\mathbb{R}^{N}}k(x) \bar{G}_{s}(u,v) \psi_{2}dx.
\end{eqnarray}
Choosing $\psi_{1}=|u^{T}|^{\nu_{1}}u^{T}$,  where  $T >\delta $, $\nu_{1}>0$, $u^{T}$ is defined as
\begin{align*}
      u^{T}
=
    \begin{cases}
         -T, \;\;\; \;\;\; \;\; \;             &\text { if }\;\;u \leq -T, \\
          u, \;\;\; \;\;\; \;\; \;             &\text { if }\;\;-T<u<T,\\
          T, \;\;\; \;\;\; \;\; \;             &\text { if }\;\;u \geq T.
    \end{cases}
\end{align*}
Choosing $\psi_{2}=|v^{T}|^{\nu_{2}}v^{T}$, where $T >\delta $, $\nu_{2}>0$, $v^{T}$ is defined by
\begin{align*}
    v^{T}
=
    \begin{cases}
         -T, \;\;\; \;\;\; \;\; \;             &\text { if }\;\;v \leq -T, \\
          v, \;\;\; \;\;\; \;\; \;             &\text { if }\;\;-T<v<T,\\
          T, \;\;\; \;\;\; \;\; \;             &\text { if }\;\;v \geq T.
    \end{cases}
\end{align*}
Since $u$ and $u^{T}$ have the same sign, by $(\Upsilon_{1})$, $(\Lambda_{0})$ and (\ref{r2.2}), it holds that
\begin{eqnarray*}
&&
         \int_{\mathbb{R}^{N}}
         \phi_{1}\left(\frac{V_{1}(x)u^{2}+|\nabla u|^{2}}{2}\right)
          (V_{1}(x)u \psi_{1}+\nabla u\cdot \nabla  \psi_{1})
          dx
\\
&=&
       \int_{\mathbb{R}^{N}}
         \phi_{1}\left(\frac{V_{1}(x)u^{2}+|\nabla u|^{2}}{2}\right)
          (V_{1}(x)u |u^{T}|^{\nu_{1}}u^{T}+\nabla u\cdot \nabla (|u^{T}|^{\nu_{1}}u^{T}))
          dx
\\
&\geq&
       \int_{\mathbb{R}^{N}}
         \phi_{1}\left(\frac{V_{1}(x)u^{2}+|\nabla u|^{2}}{2}\right)
          \nabla u\cdot \nabla (|u^{T}|^{\nu_{1}}u^{T})
          dx
\\
&\geq&
     \rho_{0} \int_{\mathbb{R}^{N}}\nabla u\cdot \nabla (|u^{T}|^{\nu_{1}}u^{T})dx
\\
&=&
 \rho_{0}(\nu_{1}+1)\int_{\mathbb{R}^{N}}|u^{T}|^{\nu_{1}}\nabla u^{T}\cdot\nabla u dx
\\
&=&
      \rho_{0}(\nu_{1}+1)
       \int_{\{x\in\mathbb{R}^{N}: -T<u(x)<T\}}
        |u^{T}|^{\nu_{1}}|\nabla u|^{2}
        dx
\\
&=&
      \rho_{0}(\nu_{1}+1)
       \int_{\{x\in\mathbb{R}^{N}: -T<u(x)<T\}}
        |u^{T}|^{\nu_{1}}|\nabla u^{T}|^{2}
        dx
\\
&=&
       \rho_{0}(\nu_{1}+1)
       \int_{\mathbb{R}^{N}}
        |u^{T}|^{\nu_{1}}|\nabla u^{T}|^{2}
        dx
\\
&=&
        \frac{4\rho_{0}(\nu_{1}+1)}{\left(\nu_{1}+2\right)^{2}}
         \int_{\mathbb{R}^{N}}
         \left|\nabla \left|u^{T}\right|^{\frac{\nu_{1}+2}{2}}\right|^{2}
         dx
\\
&\geq&
       \frac{4\rho_{0}(\nu_{1}+1)}{ D\left(\nu_{1}+2\right)^{2}}
         \left(
         \int_{\mathbb{R}^{N}}
         \left(\left|u^{T}\right|^{\frac{\nu_{1}+2}{2}}\right)^{2^{\ast}}
          dx
         \right)^{\frac{2}{2^{\ast}}},
\end{eqnarray*}
where the last step is obtained by using (\ref{r2.2}).
Then, by \eqref{3.8.1}, $(F_1')$, \eqref{A4}, H\"{o}lder's inequality, $(\Lambda_{1})$, $(F'_{5})$ and $(K)$, we have
\begin{eqnarray*}
&&
    \frac{4\rho_{0}(\nu_{1}+1)}{\left(\nu_{1}+2\right)^{2}}
         \left(
         \int_{\mathbb{R}^{N}}
         \left(\left|u^{T}\right|^{\frac{\nu_{1}+2}{2}}\right)^{2^{\ast}}
          dx
         \right)^{\frac{2}{2^{\ast}}}
\\
&\leq&
   D\int_{\mathbb{R}^{N}}  \bar{F}_{t}(x,u,v) \psi_{1}dx
    +D\varepsilon\int_{\mathbb{R}^{N}}k(x)  \bar{G}_{t}(u,v) \psi_{1}dx
\\
&\leq&
     DC_1 \int_{\mathbb{R}^{N}}|u|^{k_{1}r-1}|u^{T}|^{\nu_{1}+1}dx
    +DC_1 \int_{\mathbb{R}^{N}}|v|^{\frac{k_{2}(k_{1}r-1)}{k_{1}}}|u^{T}|^{\nu_{1}+1}dx
\\
&&
    +D|\varepsilon|\int_{\{x\in \R^{N}:|(u(x),v(x))|\leq \delta\}}
        |k(x)||\bar{G}_{t}(u,v)||u^{T}|^{\nu_{1}+1}
        dx
\\
&\leq&
     DC_1 \int_{\mathbb{R}^{N}}|u|^{k_{1}r+\nu_{1}}dx
    +DC_1 \int_{\mathbb{R}^{N}}|v|^{\frac{k_{2}(k_{1}r-1)}{k_{1}}}|u|^{\nu_{1}+1}dx
\\
&&
    +D|\varepsilon|\int_{\{x\in \R^{N}:|(u(x),v(x))|\leq \delta\}}
        |k(x)||\bar{G}_{t}(u,v)||u^{T}|^{\nu_{1}+1}
        dx
\\
&\leq&
     DC_1 \int_{\mathbb{R}^{N}}|u|^{k_{1}r+\nu_{1}}dx
    +DC_1 \int_{\mathbb{R}^{N}}
    \left(V_{2}^{-1}(x)\right)^{\frac{k_{2}(k_{1}r-1)}{2k_{1}}}
    (V_{2}(x)|v|^{2})^{\frac{k_{2}(k_{1}r-1)}{2k_{1}}}
    |u|^{\nu_{1}+1}
    dx
\\
& &
    +D|\varepsilon|\max_{|(t,s)|\leq \delta}|\bar{G}_{t}(t,s)|
      \int_{\{x\in \R^{N}:|(u(x),v(x))|\leq \delta\}}|k(x)||u|^{\nu_{1}+1}dx
\\
&\leq&
     DC_1 \int_{\mathbb{R}^{N}}|u|^{k_{1}r+\nu_{1}}dx
    +DC_1
    \left(
    \int_{\mathbb{R}^{N}}
    \left(V_{2}^{-1}(x)\right)
     ^{\frac{k_{2}(k_{1}r-1)}{2k_{1}}\cdot \frac{2k_{1}r}{(2-k_{2}r)(k_{1}r-1)}}
      dx
    \right)^{\frac{(2-k_{2}r)(k_{1}r-1)}{2k_{1}r}}
\\
&&
  \cdot\left(
    \int_{\mathbb{R}^{N}}
    (V_{2}(x)|v|^{2})^{\frac{k_{2}(k_{1}r-1)}{2k_{1}}\cdot\frac{2k_{1}}{k_{2}(k_{1}r-1)}}
     dx
    \right)^{\frac{k_{2}(k_{1}r-1)}{2k_{1}}}
    \left(
    \int_{\mathbb{R}^{N}}
    |u|^{(\nu_{1}+1)k_{1}r}
     dx
    \right)^{\frac{1}{k_{1}r} }
\\
&&
    +D\max_{|(t,s)|\leq \delta}|\bar{G}_{t}(t,s)|
      \int_{\mathbb{R}^{N}}|k(x)||u|^{\nu_{1}+1}dx
\\
&\leq&
     DC_1 \int_{\mathbb{R}^{N}}|u|^{k_{1}r+\nu_{1}}dx
    +D\|k\|_{\infty}\max_{|(t,s)|\leq \delta}|\bar{G}_{t}(t,s)|
      \int_{\mathbb{R}^{N}}|u|^{\nu_{1}+1}dx
\\
&&
    +DC_1
    \left(
    \int_{\mathbb{R}^{N}}
    \left(V_{2}^{-1}(x)\right)
     ^{\frac{k_{2}r}{2-k_{2}r}}
      dx
    \right)^{\frac{2-k_{2}r}{k_{2}r}\cdot\frac{k_{2}(k_{1}r-1)}{2k_{1}}}
    \left(
      \int_{\mathbb{R}^{N}}(|\nabla v|^{2}+V_{2}(x)|v|^{2})dx
    \right)^{\frac{k_{2}(k_{1}r-1)}{2k_{1}}}
    \left(
       \int_{\mathbb{R}^{N}}
       |u|^{(\nu_{1}+1)k_{1}r}
        dx
    \right)^{\frac{1}{k_{1}r} }
\\
&=&
     DC_1 \int_{\mathbb{R}^{N}}|u|^{k_{1}r+\nu_{1}}dx
    +DC_{6}\int_{\mathbb{R}^{N}}|u|^{\nu_{1}+1}dx
    +DC_7\|v\|_{W_{2}}^{\frac{k_{2}(k_{1}r-1)}{k_{1}}}
    \left(
       \int_{\mathbb{R}^{N}}|u|^{(\nu_{1}+1)k_{1}r}dx
    \right)^{\frac{1}{k_{1}r}},
\end{eqnarray*}
where
$$
C_{6}=\|k\|_{\infty}\max_{|(t,s)|\leq \delta}|\bar{G}_{t}(t,s)|
\;\;
\mbox{and}
\;\;
C_{7}=C_1 \|V_{2}^{-1}\|_{\frac{k_{2}r}{2-k_{2}r}}^{\frac{k_{2}(k_{1}r-1)}{2k_{1}}}.
$$
Next, we assume that
$$
\left(
\int_{\mathbb{R}^{N}}|u|^{(\nu_{1}+1)k_{1}r}dx
\right)^{\frac{1}{k_{1}r}}
=\max\left\{
\int_{\mathbb{R}^{N}}|u|^{k_{1}r+\nu_{1}}dx,
\int_{\mathbb{R}^{N}}|u|^{\nu_{1}+1}dx,
\left(
\int_{\mathbb{R}^{N}}|u|^{(\nu_{1}+1)k_{1}r}dx
\right)^{\frac{1}{k_{1}r}}
\right\}.$$
The cases
$$
\int_{\mathbb{R}^{N}}|u|^{k_{1}r+\nu_{1}}dx
=\max\left\{
\int_{\mathbb{R}^{N}}|u|^{k_{1}r+\nu_{1}}dx,
\int_{\mathbb{R}^{N}}|u|^{\nu_{1}+1}dx,
\left(
\int_{\mathbb{R}^{N}}|u|^{(\nu_{1}+1)k_{1}r}dx
\right)^{\frac{1}{k_{1}r}}
\right\}$$
and
$$
\int_{\mathbb{R}^{N}}|u|^{\nu_{1}+1}dx
=\max\left\{
\int_{\mathbb{R}^{N}}|u|^{k_{1}r+\nu_{1}}dx,
\int_{\mathbb{R}^{N}}|u|^{\nu_{1}+1}dx,
\left(
\int_{\mathbb{R}^{N}}|u|^{(\nu_{1}+1)k_{1}r}dx
\right)^{\frac{1}{k_{1}r}}
\right\}$$
can be similarly treated.
So for fixed $u\in W_{1}$, the following estimate can be obtained
$$
\frac{4\rho_{0}(\nu_{1}+1)}{\left(\nu_{1}+2\right)^{2}}
         \left(
         \int_{\mathbb{R}^{N}}
         \left(\left|u^{T}\right|^{\frac{\nu_{1}+2}{2}}\right)^{2^{\ast}}
          dx
         \right)^{\frac{2}{2^{\ast}}}
\leq
     \left(
      DC_1
      +
     DC_7\|v\|_{W_{2}}^{\frac{k_{2}(k_{1}r-1)}{k_{1}}}
      +
     DC_{6}
     \right)
    \left(
       \int_{\mathbb{R}^{N}}|u|^{(\nu_{1}+1)k_{1}r}dx
    \right)^{\frac{1}{k_{1}r}}.
$$
Note that $u^{T}\rightarrow u$ as $T\rightarrow+\infty$.
Taking the limits in estimate above, it holds
\begin{eqnarray*}
      \|u\|_{\frac{(\nu_{1}+2)N}{N-2} }
\leq
    \left(
     \frac{D_1(\nu_{1}+2)^{2}}{(\nu_{1}+1)}
     \right)^{\frac{1}{\nu_{1}+2}}
    \|u\|_{(\nu_{1}+1)k_{1}r}^{\frac{\nu_{1}+1}{\nu_{1}+2} }
\leq
    (D_1(\nu_{1}+2))^{\frac{2}{\nu_{1}+2}}
    \|u\|_{(\nu_{1}+1)k_{1}r}^{\frac{\nu_{1}+1}{\nu_{1}+2} },
\end{eqnarray*}
where
\begin{eqnarray*}
  D_1
=
     \frac{1}{4\rho_{0}}
     \left(
     DC_1
      +
     DC_7\|v\|_{W_{2}}^{\frac{k_{2}(k_{1}r-1)}{k_{1}}}
      +
     DC_{6}
     \right).
\end{eqnarray*}
Set $\nu_{1,k}=\frac{(\nu_{1,k-1}+2)N}{(N-2)k_{1}r}-1$, where $k=1,2,\cdots$ and $\nu_{1,0}=\frac{2^{\ast}-k_{1}r}{k_{1}r}$.
It is not hard to verify that
$\nu_{1,k}=\frac{\left(\frac{2^{\ast}}{2k_{1}r}\right)^{k+1}-1}{\left(\frac{2^{\ast}}{2k_{1}r}\right)-1}\nu_{1,0}$,
for $k=1,2,\cdots$, (see  \cite[Lemma A.1]{Liu-Zhang-Wang2024} for details).
Clearly, $\frac{2^{\ast}}{2k_{1}r}>1$ by $k_{1}r<\frac{2^{\ast}}{2}$.
Therefore, $\nu_{1,k}\rightarrow +\infty$ as $k\rightarrow +\infty$.
\par
By Moser's iteration method we obtain
{\footnotesize\begin{align*}
\|u\|_{\frac{(\nu_{1,k}+2)N}{N-2}}
&\leq
        \left(D_1(\nu_{1,k}+2)\right)^{\frac{2}{\nu_{1,k}+2}}
        \cdots
        \left(D_1(\nu_{1,1}+2)\right)^{\frac{2}{\nu_{1,1}+2}}
        \left(D_1(\nu_{1,0}+2)\right)^{\frac{2}{\nu_{1,0}+2}}
        \|u\|_{2^{\ast}}^{
        \frac{\nu_{1,0}+1}{\nu_{1,0}+2}
        \frac{\nu_{1,1}+1}{\nu_{1,1}+2}
        \frac{\nu_{1,2}+1}{\nu_{1,2}+2}
        \cdots
        \frac{\nu_{1,k}+1}{\nu_{1,k}+2}
        }
\nonumber\\
&
\leq
       \exp\left(\sum_{i=0}^{k}\frac{2\ln\left(D_1(\nu_{1,i}+2)\right)}{\nu_{1,i}+2}\right)
        \|u\|_{2^{\ast}}^{\zeta_{k}},
\end{align*}}
where $\zeta_{k}=\Pi_{i=0}^{k}\frac{\nu_{1,i}+1}{\nu_{1,i}+2}$.
Taking $k\rightarrow\infty$, we have
\begin{eqnarray*}
         \|u\|_{\infty}
\leq
       C^{*} \|u\|_{2^{\ast}}^{\zeta},
\end{eqnarray*}
where $0<\zeta=\Pi_{i=0}^{\infty}\frac{\nu_{1,i}+1}{\nu_{1,i}+2}<1$, $C^{*}=\exp\left(\sum_{i=0}^{\infty}\frac{2\ln\left(D_1(\nu_{1,i}+2)\right)}{\nu_{1,i}+2}\right)$ is a positive constant (for details, see \cite[Lemma A.3]{Liu-Zhang-Wang2024} and \cite[Lemma A.4]{Liu-Zhang-Wang2024}).
Similarly, by \eqref{3.8.2}, $(\Upsilon_{1})$, $(\Lambda_{0})$, $(\Lambda_{1})$, $(F_1')$, $(F'_{5})$, $(K)$, \eqref{A4} and $k_{2}r<\frac{2^{\ast}}{2}$,
we can get
\begin{eqnarray*}
         \|v\|_{\infty}
\leq
       D^{*} \|v\|_{2^{\ast}}^{\vartheta}.
\end{eqnarray*}
The proof is completed.
\qed
\par
\noindent
\begin{remark}\label{remark3.9.1}
By Lemma \ref{lemma3.8} and Remark \ref{remark2.1},
if $(u,v)$ is a critical point of $\bar{J}_{\varepsilon}(u,v)$ with
$\|u\|_{W_{1}}\leq \left(\frac{\delta}{4C^{*}S_{2^{\ast}}^{\zeta}}\right)^{\frac{1}{\zeta}}$
and
$\|v\|_{W_{2}}\leq  \left(\frac{\delta}{4D^{*}S_{2^{\ast}}^{\vartheta}}\right)^{\frac{1}{\vartheta}}$,
then
$\|u\|_{\infty}\leq \frac{\delta}{4}$
and
$\|v\|_{\infty}\leq \frac{\delta}{4}$.
This means that $(u,v)$ is solution of the original problem  \eqref{eq1}.
\end{remark}
\par
\noindent
\begin{lemma}\label{lemma3.7}
Assume that  $(\Upsilon_{1})$--$(\Upsilon_{3})$, $(F'_{1})$, $(F'_{2})$, $(F'_{5})$, $(\Lambda_{0})$ and $( K)$ hold.
For any $b>0$, there exist  $\sigma(b)>0$ such that if $\bar{J}_{\varepsilon}'(u,v)=0$ and $|\bar{J}_{\varepsilon}(u,v)|\leq \sigma(b)$, where $|\varepsilon|\leq \sigma(b)$,
then $\|(u,v)\|\leq b$.
\end{lemma}
\par
\noindent
{\bf Proof.}
Suppose on the contrary that there exist two sequences $\{(u_{n},v_{n})\}\subset W$ and $\{\varepsilon_{n}\}$ such that $\varepsilon_{n}\rightarrow 0$,
$\bar{J}_{\varepsilon_{n}}(u_{n},v_{n})\rightarrow 0$ as $n\rightarrow\infty$, $\bar{J}_{\varepsilon_{n}}'(u_{n},v_{n})=0$ and $\|(u_{n},v_{n})\|\geq b_{0}>0$, where $b_{0}$ is independent of $n$.
Then $\{(u_{n},v_{n})\}$ is the $(PS)_{0}$ sequence of $\bar{J}_{0}$.
Thus it follows from Lemma \ref{lemma3.5} that a subsequence of $\{(u_{n},v_{n})\}$ which converges to $(u_{0},v_{0})$ in $W$, satisfies
\begin{eqnarray*}
0=     \langle \bar{J}_{0}'(u_{0},v_{0}),(u_{0},v_{0})\rangle
& =&
       \int_{\mathbb{R}^{N}}\phi_{1}
        \left(\frac{V_{1}(x)u_{0}^{2}+|\nabla u_{0}|^{2}}{2}\right)
        (V_{1}(x)u_{0}^{2}+|\nabla u_{0}|^{2})dx
\nonumber\\
 & &
       +\int_{\mathbb{R}^{N}}\phi_{2}
        \left(\frac{V_{2}(x)v_{0}^{2}+|\nabla v_{0}|^{2}}{2}\right)
        (V_{2}(x)v_{0}^{2}+|\nabla v_{0}|^{2})dx
\nonumber\\
&    &
       -\int_{\mathbb{R}^{N}} \bar{F}_{t}(x,u_{0},v_{0})u_{0}dx
       - \int_{\mathbb{R}^{N}}\bar{F}_{s}(x,u_{0},v_{0})v_{0}dx\nonumber\\
\end{eqnarray*}
and
{\small\begin{eqnarray*}
     \bar{J}_{0}(u_{0},v_{0})
=       \int_{\mathbb{R}^{N}}\left\{\Phi_{1}\left(\frac{V_{1}(x)u_{0}^{2}+|\nabla u_{0}|^{2}}{2}\right)
        +\Phi_{2}\left(\frac{V_{2}(x)v_{0}^{2}+|\nabla v_{0}|^{2}}{2}\right)\right\}dx
        - \int_{\mathbb{R}^{N}}\bar{F}(x,u_{0},v_{0})dx
=0.
\end{eqnarray*}}
From the above two equations, $(\Upsilon_{3})$, $(\Upsilon_{1})$, $(F'_{2})$ and $\beta<2$, it holds
\begin{eqnarray*}\label{3.7.1}
&&
0=     \bar{J}_{0}(u_{0},v_{0})-\frac{1}{\beta}\langle \bar{J}_{0}'(u_{0},v_{0}),(u_{0},v_{0})\rangle
\nonumber\\
&&\quad
=
       \int_{\mathbb{R}^{N}}
       \left\{
       \Phi_{1}\left(\frac{V_{1}(x)u_{0}^{2}+|\nabla u_{0}|^{2}}{2}\right)
        -\frac{2}{\beta}\phi_{1}\left(\frac{V_{1}(x)u_{0}^{2}+|\nabla u_{0}|^{2}}{2}\right)
        \frac{V_{1}(x)u_{0}^{2}+|\nabla u_{0}|^{2}}{2}
        \right\}
        dx
\nonumber\\
&&\qquad
    +\int_{\mathbb{R}^{N}}
       \left\{
       \Phi_{2}\left(\frac{ V_{2}(x)v_{0}^{2}+|\nabla v_{0}|^{2}}{2}\right)
        -\frac{2}{\beta}\phi_{2}\left(\frac{V_{2}(x)v_{0}^{2}+|\nabla v_{0}|^{2}}{2}\right)
        \frac{V_{2}(x)v_{0}^{2}+|\nabla v_{0}|^{2}}{2}
        \right\}
        dx
\nonumber\\
&&\qquad
     +\int_{\mathbb{R}^{N}}
     \left(
     \frac{1}{\beta}\bar{F}_{t}(x,u_{0},v_{0})u_{0}
     +\frac{1}{\beta}\bar{F}_{s}(x,u_{0},v_{0})v_{0}
     -\bar{F}(x,u_{0},v_{0})
     \right)dx
\nonumber\\
&&\quad
\leq
      \frac{\beta-2}{\beta}
      \int_{\mathbb{R}^{N}}\Phi_{1}\left(\frac{V_{1}(x)u_{0}^{2}+|\nabla u_{0}|^{2}}{2}\right)dx
      +
      \frac{\beta-2}{\beta}
      \int_{\mathbb{R}^{N}}\Phi_{2}\left(\frac{V_{2}(x)v_{0}^{2}+|\nabla v_{0}|^{2}}{2}\right)dx
\nonumber\\
&&\qquad
     +\int_{\mathbb{R}^{N}}
     \left(
     \frac{1}{\beta}\bar{F}_{t}(x,u_{0},v_{0})u_{0}
     +\frac{1}{\beta}\bar{F}_{s}(x,u_{0},v_{0})v_{0}
     -\bar{F}(x,u_{0},v_{0})
     \right)dx
\nonumber\\
&&\quad
\leq
      \frac{\rho_{0}(\beta-2)}{2\beta}
      \int_{\mathbb{R}^{N}}\left(V_{1}(x)u_{0}^{2}+|\nabla u_{0}|^{2}\right)dx
      +
      \frac{\rho_{0}(\beta-2)}{2\beta}
      \int_{\mathbb{R}^{N}}\left(V_{2}(x)v_{0}^{2}+|\nabla v_{0}|^{2}\right)dx
\nonumber\\
&&\qquad
     +\int_{\mathbb{R}^{N}}
     \left(
     \frac{1}{\beta}\bar{F}_{t}(x,u_{0},v_{0})u_{0}
     +\frac{1}{\beta}\bar{F}_{s}(x,u_{0},v_{0})v_{0}
     -\bar{F}(x,u_{0},v_{0})
     \right)dx
\nonumber\\
&&\quad
\leq
      \frac{\rho_{0}(\beta-2)}{2\beta}
      \int_{\mathbb{R}^{N}}\left(V_{1}(x)u_{0}^{2}+V_{2}(x)v_{0}^{2}\right)dx
     +\int_{\mathbb{R}^{N}}
       \left(
        \frac{1}{\beta}\bar{F}_{t}(x,u_{0},v_{0})u_{0}
         +\frac{1}{\beta}\bar{F}_{s}(x,u_{0},v_{0})v_{0}
         -\bar{F}(x,u_{0},v_{0})
       \right)dx
\nonumber\\
&&\quad
\leq
      \frac{\rho_{0}(\beta-2)}{2\beta}
      \int_{\mathbb{R}^{N}}\left(V_{1}(x)u_{0}^{2}+V_{2}(x)v_{0}^{2}\right)dx
     +
     \frac{\rho_{0}(2-\beta)}{4\beta}
     \int_{\mathbb{R}^{N}}\left(V_{1}(x)u_{0}^{2}+V_{2}(x)v_{0}^{2}\right)dx
\nonumber\\
&&\quad
=-\frac{\rho_{0}(2-\beta)}{4\beta}
     \int_{\mathbb{R}^{N}}\left(V_{1}(x)u_{0}^{2}+V_{2}(x)v_{0}^{2}\right)dx,
\end{eqnarray*}
which implies $(u_{0},v_{0})\equiv(0,0)$ due to the fact that $\beta\in(0,2)$ and $(\Lambda_{0})$.
However, from $\|(u_{n},v_{n})\|\geq b_{0}>0$ and $(u_{n},v_{n})\rightarrow (u_{0},v_{0})$ in $W$, we have
$\|(u_{0},v_{0})\|\geq b_{0}>0$, which contradicts the fact that $(u_{0},v_{0})\equiv (0,0)$. The proof is completed.
\qed

\par
\noindent
\begin{remark}\label{remark3.7.1}
By Lemma \ref{lemma3.7}, if $(u,v)$ is a critical point of $\bar{J}_{\varepsilon}(u,v)$ with
$|\varepsilon|
\leq
\sigma\left(
\min\left\{
\left(\frac{\delta}{4C^{*}S_{2^{\ast}}^{\zeta}}\right)^{\frac{1}{\zeta}},
\left(\frac{\delta}{4D^{*}S_{2^{\ast}}^{\vartheta}}\right)^{\frac{1}{\vartheta}}
\right\}
\right)$
and
$|\bar{J}_{\varepsilon}(u,v)|
\leq
\sigma\left(
\min\left\{
\left(\frac{\delta}{4C^{*}S_{2^{\ast}}^{\zeta}}\right)^{\frac{1}{\zeta}},
\left(\frac{\delta}{4D^{*}S_{2^{\ast}}^{\vartheta}}\right)^{\frac{1}{\vartheta}}
\right\}
\right)$,
then
$\|(u,v)\|
\leq
\min\left\{
\left(\frac{\delta}{4C^{*}S_{2^{\ast}}^{\zeta}}\right)^{\frac{1}{\zeta}},
\left(\frac{\delta}{4D^{*}S_{2^{\ast}}^{\vartheta}}\right)^{\frac{1}{\vartheta}}
\right\}.
$
\end{remark}

\par
\vskip2mm
\noindent
{\bf  Proof of Theorem \ref{theorem1.1}}
Without loss of generality, we assume that $\varepsilon>0$.
The case $\varepsilon<0$ can be studied similarly by replacing $\bar{G}_u(u,v)$ and $\bar{G}_v(u,v)$
with $-\bar{G}_u(u,v)$ and $-\bar{G}_v(u,v)$.
Next, we are ready to verify that $\bar{J}_{\varepsilon}(u,v)$ satisfies conditions $(A_{1})$, $(A_{2})$, $(A_{4})$ and  $(A_{5})$ in Lemma \ref{lemma2.2} and condition $(A_{3})'$ in Remark \ref{remark2.2}.
To verify condition $(A_{1})$, by \eqref{3.4.2}, we have
\begin{eqnarray*}
 \inf_{\varepsilon\in [0,1],(u,v)\in W}\bar{J}_{\varepsilon}(u,v)>-\infty.
\end{eqnarray*}
By $(F'_{5})$ and $(K)$, we have
\begin{eqnarray}\label{3.4.1}
      |\bar{J}_{\varepsilon}(u,v)-\bar{J}_{0}(u,v)|
&\leq&
      |\varepsilon|\int_{\mathbb{R}^{N}}|k(x)||\widetilde{G}(u,v)|dx
\nonumber\\
&=&
        |\varepsilon|\int_{\{x\in \R^{N}:|(u(x),v(x))|\leq \delta\}}|k(x)||\bar{G}(u,v)|dx
\nonumber\\
&\leq&
       |\varepsilon|\max_{|(t,s)|\leq \delta}|\bar{G}(t,s)|\int_{|(u,v)|\leq \delta}|k(x)|dx
\nonumber\\
&\leq&
       |\varepsilon|\max_{|(t,s)|\leq \delta}|\bar{G}(t,s)|\int_{\mathbb{R}^{N}}|k(x)|dx
\nonumber\\
&\leq&
     |\varepsilon|C:=\psi(\varepsilon)
\end{eqnarray}
where $C$ is a constant independent of $(u,v)$ and $\varepsilon$.
So, condition $(A_{2})$ holds.
Condition $(A_{4})$ holds since $(F'_{0})$ and $(F'_{4})$.
Conditions $(A_{3})'$ and $(A_{5})$ follow from Lemma \ref{lemma3.5} and Lemma \ref{lemma3.6}, respectively.
Thus, $\bar{J}_{\varepsilon}(u,v)$ satisfies all the conditions in Lemma \ref{lemma2.2}.
By arguments as in \cite[Proof of Theorem 1.1]{Kajikiya2013} or \cite[Proof of Corollary 1.1]{Huang2022}, for any $\sigma>0$ and any given $k\in \mathbb N$, we have $k$ distinct critical values of $\bar{J}_{\varepsilon}$ satisfying
\begin{eqnarray}\label{z.1}
-\sigma
<   a_{n(1)}(\varepsilon)
<   a_{n(2)}(\varepsilon)
<
\cdots
<   a_{n(k)}(\varepsilon)
<0.
\end{eqnarray}
Finally, due to the arbitrariness of $\sigma$,
we can take
$$
0
<
\sigma
\leq
\sigma\left(
\min\left\{
\left(\frac{\delta}{4C^{*}S_{2^{\ast}}^{\zeta}}\right)^{\frac{1}{\zeta}},
\left(\frac{\delta}{4D^{*}S_{2^{\ast}}^{\vartheta}}\right)^{\frac{1}{\vartheta}}
\right\}
\right).$$
Then by Remark \ref{remark3.7.1} and Remark \ref{remark3.9.1}, the original problem  \eqref{eq1} has at least $k$ solutions whose $L^{\infty}$-norms are less than $\frac{\delta}{2}$. \qed

\vskip2mm
\section{Example}\label{section 4}
Consider this elliptic system
\begin{eqnarray}\label{ex1}
 \begin{cases}
 -\mbox{div}\left\{\left(\left(1+\frac{V_{1}(x)u^{2}+|\nabla u|^{2}}{2}\right)^{-2}+4\right)\nabla u\right\}
 +\left(\left(1+\frac{V_{1}(x)u^{2}+|\nabla u|^{2}}{2}\right)^{-2}+4\right)V_{1}(x)u
 = F_u(x,u,v), \ \ x\in \mathbb{R}^6,\\
 -\mbox{div}\left\{\left(\left(1+\frac{V_{2}(x)v^{2}+|\nabla v|^{2}}{2}\right)^{-3}+6\right)\nabla v\right\}
 +\left(\left(1+\frac{V_{2}(x)v^{2}+|\nabla v|^{2}}{2}\right)^{-3}+6\right)V_{2}(x)v
= F_v(x,u,v), \ \ x\in \mathbb{R}^6,\\
 \end{cases}
\end{eqnarray}
where
\begin{eqnarray}\label{4.1.1}
F(x,t,s)= (\sin x_{1}+3)(|t|^{\frac{3}{2}}+|s|^{\frac{3}{2}}+|t|^{2}|s|^{2}).
\end{eqnarray}
\par
Let $\phi_1(s)=(1+s)^{-2}+4$ and $\phi_2(s)=(1+s)^{-3}+6$ for $s\geq0$.
Then $\Phi_1(s)=-(1+s)^{-1}+4s$, $\Phi_2(s)=-\frac{1}{2}(1+s)^{-2}+6s$ for $s\geq0$.
By simple computations we have that $\phi_i (i=1,2)$ satisfy $(\Upsilon_1)$--$(\Upsilon_3)$, $\rho_{0}=4$ and  $\rho_{1}=7$.
\par
 Let $V_1(x)=\sum_{i=1}^{6}x_i^{2}+20$ and $V_2(x)=\sum_{i=1}^{6}x_i^{4}+30$ for all $x\in \mathbb{R}^{6}$.
It is easy to verify that $V_i (i=1,2)$ satisfy $(\Lambda_{0})$ and $(\Lambda_{1})$.
 \par
It follows from (\ref{4.1.1}) that
\begin{eqnarray*}
&F_{t}(x,t,s)=(\sin x_{1}+3)\left(\frac{3}{2}t|t|^{-\frac{1}{2}}+2t|s|^{2}\right), \label{4.1.2}\\
&F_{s}(x,t,s)=(\sin x_{1}+3)\left(\frac{3}{2}s|s|^{-\frac{1}{2}}+2s|t|^{2}\right).\label{4.1.3}
\end{eqnarray*}
Thus,
\begin{eqnarray*}
&|F_{t}(x,t,s)|\leq 8(t^{\frac{1}{4}}+ s^{\frac{1}{4}}), \;\;\text{for\;all}\;\;|(t,s)|\leq1, \label{4.1.4}\\
&|F_{s}(x,t,s)|\leq 8(t^{\frac{1}{4}}+ s^{\frac{1}{4}}), \;\;\text{for\;all}\;\;|(t,s)|\leq1,\label{4.1.5}
\end{eqnarray*}
where we choose $\delta=1$, $C_{1}=C_{2}=8$, $k_{1}=k_{2}=2<2^{\ast}=3$ and $r=\frac{5}{8}$.
Furthermore, by Young's inequality, it holds
\begin{align*}
&
tF_{t}(x,t,s)+sF_{s}(x,t,s)-\beta F(x,t,s)\\
&
=
(\sin x_{1}+3)\left(\frac{3}{2}-\beta\right)\left(|t|^{\frac{3}{2}}+|s|^{\frac{3}{2}}\right)
+ (\sin x_{1}+3)(4-\beta)|t|^{2}|s|^{2}
     \\
&
\leq   9|t|^{2}|s|^{2}
\leq   \frac{9}{2}t^{4}+\frac{9}{2}s^{4}
\leq      \frac{9}{2}|t|^{2}+\frac{9}{2}|s|^{2}
\\
&
\leq
\frac{\rho_{0}(2-\beta)}{2}(V_{1}(x)|t|^{2}+V_{2}(x)|s|^{2}), \;\;\text{for\;all}\;\;|(t,s)|\leq1,
\end{align*}
where $\beta=\frac{7}{4}<2$.
We easily see that conditions $(F_3)$ and $(F_4)$ hold.
Since
\begin{align*}
      \lim_{|(t,s)|\rightarrow 0}
      \left(\inf_{x\in \mathbb{R}^{N}}\frac{F(x,t,s)}{|t|^{\frac{7}{4}}+|s|^{\frac{7}{4}}}\right)
&=     \lim_{|(t,s)|\rightarrow 0}
      \left(\inf_{x\in \mathbb{R}^{N}}
      \frac{(\sin x_{1}+3)\left(|t|^{\frac{3}{2}}+|s|^{\frac{3}{2}}+t^{2}s^{2}\right)}
      {|t|^{\frac{7}{4}}+|s|^{\frac{7}{4}}}
      \right)
=     \lim_{|(t,s)|\rightarrow 0}
       \frac{2\left(|t|^{\frac{3}{2}}+|s|^{\frac{3}{2}}+t^{2}s^{2}\right)}
      {|t|^{\frac{7}{4}}+|s|^{\frac{7}{4}}}\\
&\leq  \lim_{|(t,s)|\rightarrow 0}2\left(|t|^{-\frac{1}{4}}+|s|^{-\frac{1}{4}}\right)
      +\lim_{|(t,s)|\rightarrow 0}\frac{2t^{2}s^{2}}{|t|^{\frac{7}{4}}+|s|^{\frac{7}{4}}}
=+\infty.
\end{align*}
So the condition $(F_4)$ holds. Then by Corollary \ref{corollary1.1}, system (\ref{ex1}) has infinitely many nontrivial weak solutions.
\section{Results for the scalar equation}\label{section 5}
In this section, we consider the multiplicity of solutions for the following non-homogeneous elliptic equation with perturbation:
\begin{equation}\label{s5-eq1}
 \left\{
  \begin{array}{ll}
 -\mbox{div}\left\{\phi \left(\frac{V(x)u^{2}+|\nabla u|^{2}}{2}\right)\nabla u\right\}
                   +\phi \left(\frac{V(x)u^{2}+|\nabla u|^{2}}{2}\right)V(x)u
 =f(x,u)+\varepsilon k(x)g(u),\;\; x\in \mathbb{R}^{N},
 \\
 u\in H^{1}(\mathbb{R}^{N}),
  \end{array}
 \right.
 \end{equation}
where $N>2$ is an integer, $\phi:[0,+\infty)\rightarrow \mathbb{R}$ is a continuous function which satisfies the following conditions:
\begin{itemize}
\item[$(\Upsilon_{1})'$] there exist two constants $0< \rho_0\leq \rho_1$ such that
 $ 0<\rho_0\leq \phi(s)\leq \rho_1$ for all $s\in [0,+\infty)$;
\item[$(\Upsilon_{2})'$] let $\hbar(s):=\Phi(s^{2})$ where $\Phi(s):=\int_{0}^{s}\phi(\varsigma)d\varsigma$, there exists $l>0$ such that
$\hbar(t)\geq \hbar(s)+\hbar'(s)(t-s)+l(t-s)^{2}$
for all $t,s\geq0$;
\item[$(\Upsilon_{3})'$] $\phi(s)s\geq \Phi(s)$ for all $s\in [0,+\infty)$.
\end{itemize}
Furthermore, to state our result, we introduce the following  assumptions about $f$, $g$, $V$ and $k$:
\begin{itemize}
\item[$(f_0)$]
$f: \mathbb{R}^N\times [-\delta, \delta]\rightarrow \mathbb{R}$ is a continuous function for some $\delta>0$;
\item[$(f_{1})$] there exists a constant $c_1>0$ such that
\begin{equation*}
|f(x,t)|\leq c_1|t|^{ar-1}
 \end{equation*}
for all $|t|\leq\delta$  and $x\in  \mathbb{R}^{N}$, where $2\leq a<2^{\ast}$, $r$ is a constant with $0<r<1$ such that $1<ar<\min\{2,\frac{2^{\ast}}{2}\}$;
\item[$(f_{2})$] there exists a constant $\beta\in(0,2)$ such that
$$
tf(x,t)-\beta F(x,t)<\frac{\rho_{0}(2-\beta)}{4}V(x)|t|^{2}
$$
for all $|t|\leq\delta$  and $x\in  \mathbb{R}^{N}$;
\item[$(f_{3})$]
$$\lim_{|t|\rightarrow0}
\left(\inf_{x\in \mathbb{R}^{N}}\frac{F(x,t)}{|t|^{\beta}} \right)=+\infty;$$
\item[$(f_4)$]
 $f(x,-t)=-f(x,t)$
for all $|t|\leq\delta$  and $x\in  \mathbb{R}^{N}$.
\end{itemize}
\begin{itemize}
\item[$(g)$]
$g: [-\delta, \delta] \rightarrow \mathbb{R}$ is a continuous function.
\end{itemize}
\begin{itemize}
\item[$(\Lambda_{0})'$]
$V\in C(\mathbb{R}^{N},\mathbb{R})$
and there exists a constant $V_{0}$ such that
$\inf_{x \in \mathbb{R}^{N} }V(x)=V_{0}>0$;
\item[$(\Lambda_{1})'$]
the function $[V(x)]^{-1}$ belongs to $L^{\frac{ar}{2-ar}}(\mathbb{R}^{N})\cap L^{1}(\mathbb{R}^{N})$.
\end{itemize}
\begin{itemize}
\item[$(K)'$]
$k(x)\in L^{1}(\mathbb{R}^N)\cap L^{\infty}(\mathbb{R}^N)$.
\end{itemize}
\par
By similar proofs of Theorem \ref{theorem1.1}, we can obtain the following result.
\par
\noindent
\begin{theorem}\label{theorem5.1}
Assume that $(\Upsilon_{1})'$--$(\Upsilon_{3})'$, $(f_0)$--$(f_4)$, $(g)$, $(\Lambda_{0})'$, $(\Lambda_{1})'$ and $(K)'$ hold.
Then for any $k\in\mathbb{N}$, there exists an constant $\varepsilon(k)>0$ such that if $|\varepsilon|\leq \varepsilon(k)$, equation (\ref{s5-eq1}) possesses at least $k$ distinct solutions whose $L^{\infty}$-norms are less than $\frac{\delta}{2}$.
\end{theorem}
\par
By Theorem \ref{theorem5.1}, it is easy to obtain the following corollary.
\par
\noindent
\begin{corollary}\label{corollary5.1}
 Assume that   $(\Upsilon_{1})'$--$(\Upsilon_{3})'$, $(f_0)$--$(f_4)$, $(\Lambda_{0})'$ and $(\Lambda_{1})'$ hold.
 Then equation (\ref{s5-eq1}) with $\varepsilon=0$ possesses infinitely many distinct solutions whose $L^{\infty}$-norms are less than $\frac{\delta}{2}$.
\end{corollary}

 \vskip2mm
 \noindent
 {\bf Acknowledgments}\\
This work is supported by Yunnan Fundamental Research Projects of China (grant No: 202301AT070465) and  Xingdian Talent
Support Program for Young Talents of Yunnan Province in China.

\vskip2mm
\renewcommand\refname{References}
{}

\end{document}